\documentclass[a4paper,10p]{amsart}
\usepackage[foot]{amsaddr}
\usepackage{geometry}
\theoremstyle{remark}

\usepackage{graphicx}
\usepackage{lineno}
\usepackage{bm}
\usepackage{amsfonts} 
\usepackage{amsmath}
\usepackage{amssymb}
\usepackage{mathrsfs}
\usepackage{framed} % Framing content
\usepackage{multicol} % Multiple columns environment
\usepackage{tabularx}
\usepackage{tabulary}
\usepackage{siunitx}
\usepackage{caption}
\usepackage{float}
\usepackage{url} 
\usepackage{multicol}
\usepackage{breakurl} 

\usepackage{subcaption}
\usepackage{xpatch}
\usepackage{nomencl} % Nomenclature package
\makenomenclature
\makenomenclature 
\RequirePackage{ifthen}
\let\oldref\ref
\renewcommand{\ref}[1]{(\oldref{#1})}

 \graphicspath{
    {figures/}   
}

\usepackage{amssymb}

\newcommand{\ltwonorm}[1]{\langle #1 \rangle_{L_2(\Omega)}}

\DeclareMathAlphabet\mathbfcal{OMS}{cmsy}{b}{n}

\xpatchcmd{\thenomenclature}{%
  \section*{\nomname}% Look for `\section*... etc.
}{% Replace it by 'nothing'
}{\typeout{Success}}{\typeout{Failure}}

\usepackage{ifthen}
\renewcommand{\nomgroup}[1]{%
  \ifthenelse{\equal{#1}{A}}{\item[\textbf{Abbreviations}]}{%
    \ifthenelse{\equal{#1}{G}}{\item[\textbf{Symbols}]}{%
      \ifthenelse{\equal{#1}{C}}{\item[\textbf{Abbreviations}]}{%
        \ifthenelse{\equal{#1}{S}}{\item[\textbf{Subscripts}]}{%
          \ifthenelse{\equal{#1}{Z}}{\item[\textbf{Mathematical Symbols}]}{}
        }% matches mathematical symbols
      }% matches Subscripts
    }% matches Abbreviations
  }% matches Greek Symbols
}% matches Roman Symbols

\begin{document}
% Abbreviations
\nomenclature{$\text{ROM}$}{Reduced Order Model}
\nomenclature{$\text{POD}$}{Proper Orthogonal Decomposition}
\nomenclature{$\text{PODI}$}{Proper Orthogonal Decomposition with Interpolation}
\nomenclature{$\text{RBF}$}{Radial Basis Functions}

% Symbols
\nomenclature[G]{$\bm{u}$}{velocity field}
\nomenclature[G]{$\bm{u}'$}{fluctuating velocity field}
\nomenclature[G]{$\bar{\bm{u}}$}{mean velocity field}
\nomenclature[G]{$p'$}{fluctuating pressurefield}
\nomenclature[G]{$\bar{p}$}{mean pressure field}
\nomenclature[G]{$\theta'$}{fluctuating temperature field}
\nomenclature[G]{$\bar{\theta}$}{mean temperature field}
\nomenclature[G]{$\bm{u}'$}{fluctuating velocity field}
\nomenclature[G]{$\bar{\bm{u}}$}{mean velocity field}
\nomenclature[G]{$\bm{f}$}{Dirichlet boundary condition for velocity}
\nomenclature[G]{$g$}{Dirichlet boundary condition for temperature}
\nomenclature[G]{$\bm{h}$}{initial condition for velocity}
\nomenclature[G]{$e$}{initial condition for temperature}
\nomenclature[G]{${p}$}{pressure field}
\nomenclature[G]{${\theta}$}{temperature field}
\nomenclature[G]{${\nu}$}{dimensionless kinematic viscosity}
\nomenclature[G]{${\nu_t}$}{eddy viscosity}
\nomenclature[G]{${Pr_t}$}{turbulent Prandtl number}
\nomenclature[G]{${\alpha_{dif}}$}{thermal diffusivity}
\nomenclature[G]{${\alpha_{dif_t}}$}{turbulent thermal diffusivity}
\nomenclature[G]{${N_u^h}$}{number of degrees of freedom for velocity at full-order level}
\nomenclature[G]{${N_u^s}$}{number of uknowns for velocity at reduced level before the truncation}
\nomenclature[G]{${N_p^s}$}{number of unknowns for pressure at reduced order level before the truncation}
\nomenclature[G]{${N_{\theta}^s}$}{number of unknowns for temperature at reduced order level before the truncation}
\nomenclature[G]{${N_u^r}$}{number of unknowns for velocity at reduced order level after the truncation}
\nomenclature[G]{$C$}{correlation matrix}
\nomenclature[G]{$\bm{W}$}{eigenvector matrix}
\nomenclature[G]{$N_t$}{number of time instances}
\nomenclature[G]{$Q$}{space-time domain}
\nomenclature[G]{$T$}{final time}
% \nomenclature[G]{$\mathcal{T}$}{tessellation}
\nomenclature[G]{$\mathcal{P}$}{parameter space}
\nomenclature[G]{$\mathcal{K}$}{training set space}
\nomenclature[G]{${\Omega}$}{bounded domain}
\nomenclature[G]{${\Gamma}$}{boundary of $\Omega$}
\nomenclature[G]{${\epsilon}_{L^2}$}{$L^2$ norm error}
\nomenclature[G]{$\bm{n}$}{outward normal vector}
\nomenclature[G]{$\bm{\varphi_i}$}{i-th POD basis function for velocity}

\nomenclature[G]{${\psi_i}$}{i-th POD basis function for pressure}
\nomenclature[G]{${\chi_i}$}{i-th POD basis function for temperature}
\nomenclature[G]{${\xi_i}$}{i-th POD basis function for eddy viscosity}
% \nomenclature[G]{$L_u$}{reduced basis space for velocity}
% \nomenclature[G]{$L_p$}{reduced basis space for pressure}
% \nomenclature[G]{$L_t$}{reduced basis space for temperature}
% \nomenclature[G]{$L_s$}{reduced basis space for supremizers}
\nomenclature[G]{$\bm{M}$}{ROM mass matrix} 
\nomenclature[G]{$\bm{Q}$}{ROM convection tensor}
\nomenclature[G]{$\bm{Q}_{T1}$}{ROM turbulent tensor}
\nomenclature[G]{$\bm{Q}_{T2}$}{ROM turbulent tensor}
\nomenclature[G]{$\bm{P}$}{ROM pressure gradient matrix}
\nomenclature[G]{$\bm{K}$}{ROM mass matrix for the heat equation}
\nomenclature[G]{$\bm{G}$}{ROM convection matrix for the heat equation}
\nomenclature[G]{$\bm{N}$}{ROM diffusion matrix for the heat equation}   
\nomenclature[G]{${N_{\mu}}$}{number of parameters in the training set $\mathcal{K}$}
\nomenclature[G]{$\otimes$}{tensor product}
\nomenclature[G]{$\bm{\nabla}\cdot$}{divergence operator}
\nomenclature[G]{$\bm{\nabla}\times$}{curl operator}
\nomenclature[G]{$\bm{\nabla}$}{gradient operator}
\nomenclature[G]{$\bm{\nabla}^s$}{symmetric gradient operator}
\nomenclature[G]{$\bm{\alpha}$}{reduced vector of unknowns for velocity}
\nomenclature[G]{$\bm{b}$}{reduced vector of unknowns for pressure}
\nomenclature[G]{$\bm{c}$}{reduced vector of unknowns for temperature}
\nomenclature[G]{$\bm{l}$}{reduced vector of unknowns for eddy viscosity}
\nomenclature[G]{$\Delta$}{laplacian operator}
\nomenclature[G]{$\left\lVert \cdot\right\rVert$}{norm in $L^2(\Omega)$}
\nomenclature[G]{$\langle \cdot , \cdot \rangle$}{inner product in $L^2(\Omega)$}
\nomenclature[G]{$\tau^R$}{Reynolds stress tensor}
\nomenclature[G]{$H^R$}{Heat flux term}
\nomenclature[G]{$\bm{\overline{D}}^R$}{Reynolds-averaged strain rate tensor}
\nomenclature[G]{$\bm{I}$}{Identity matrix}
\nomenclature[G]{$k$}{Turbulent kinetic energy}
\nomenclature[G]{$\omega$}{Specific dissipation rate}
\nomenclature[G]{$\epsilon$}{Dissipation rate}
\nomenclature[G]{$\bm{U}_s$}{Snapshot matrix}
\nomenclature[G]{$\bm{w}$}{Vector of weights}
\nomenclature[G]{$\bm{\Theta}$}{Radial Basis Function kernels}
\nomenclature[G]{$\gamma$}{Spread of a kernel}
\nomenclature[G]{$\zeta$}{Control functions}
\nomenclature[G]{$u_{D}$}{Scaling coefficients}
\nomenclature[G]{$\bm{U}_{sn}$}{Global snapshot matrix resulting from nested POD}
\nomenclature[G]{$\bm{U}_{\mbox{ nested }}^i$}{i-th nested snapshot matrix}
\nomenclature[G]{$\epsilon_{L^2}$}{Relative error}

\newcommand{\dif}{\mbox{d}}

\title{A Hybrid Reduced Order Method for Modelling \protect\\ Turbulent Heat Transfer Problems\protect\\ } 

% \author[sis]{Giovanni Stabile\corref{cor1}}
% \ead{gstabile@sissa.it} 
% \author[sis]{Gianluigi Rozza} 
% \ead{grozza@sissa.it} 

% \address[sis]{SISSA, International School for Advanced Studies, Mathematics Area, mathLab, via Bonomea 265 – 34136, Trieste, Italy}

\author{Sokratia Georgaka\textsuperscript{1,*}}
\thanks{\textsuperscript{*}Corresponding Author.}
\email{s.georgaka16@imperial.ac.uk}

\author{Giovanni Stabile\textsuperscript{2}}
\email{gstabile@sissa.it}

\author{Kelbij Star\textsuperscript{3,}\textsuperscript{4}}
\email{kelbij.star@sckcen.be}

\author{Gianluigi Rozza\textsuperscript{2}}
\email{grozza@sissa.it}

\author{Michael J. Bluck\textsuperscript{1}}
\email{m.bluck@imperial.ac.uk}

\address{\textsuperscript{1}Imperial College London, Department of Mechanical Engineering, London, SW7 2BX, UK.}
\address{\textsuperscript{2}SISSA, International School for Advanced Studies, Mathematics Area, mathLab Trieste, Italy.}
\address{\textsuperscript{3}Institute for Advanced Nuclear Systems, SCK$\cdot$CEN, Mol, Belgium.}
\address{\textsuperscript{4}Department of Flow, Heat and Combustion Mechanics, Ghent University, Ghent, Belgium.}



\date{}

\dedicatory{}

\begin{abstract}
A parametric, hybrid reduced order model approach based on the Proper Orthogonal Decomposition with both Galerkin projection and interpolation based on Radial Basis Functions method is presented. This method is tested against a case of turbulent non-isothermal mixing in a T-junction pipe, a common flow arrangement found in nuclear reactor cooling systems. The reduced order model is derived from the 3D unsteady, incompressible Navier-Stokes equations weakly coupled with the energy equation. For high Reynolds numbers, the eddy viscosity and eddy diffusivity are incorporated into the reduced order model with a Proper Orthogonal Decomposition (nested and standard) with Interpolation (PODI), where the interpolation is performed using Radial Basis Functions. The reduced order solver, obtained using a $k$-$\omega$ SST URANS full order model, is tested against the full order solver in a 3D T-junction pipe with parametric velocity inlet boundary conditions.
\end{abstract}
\maketitle

\begin{multicols}{2}

%%      INTRODUCTION - sec:intro
%
%%%%%%%%%%%%%%%%%%%%%%%%%%%%%%%%%
\section{Introduction}\label{sec:intro}
TThe simulation of realistic large-scale many-query systems and the use of real-time control is an ambition in almost every industry. These large-scale systems are usually governed by partial differential equations (PDEs) which are complex and generally nonlinear. In the case of parametric PDEs, which are essential in design, control and optimisation, such large scale systems need to be solved for a large range of parameter values. As a consequence, their solution requires great computational effort, even with the use of state-of-the art computers and sophisticated computational libraries. A motivating example of this type of challenge can be found in the nuclear industry where a variety of phenomena, including turbulent flows, multiphase flows and heat transfer occur in complex geometries and therefore, obtaining high fidelity solutions across a wide range of scenarios becomes an unfeasible task without considerable simplification.        

Reduced Order Models (ROMs) have been proposed as a way of overcoming the computational burden required to obtain high fidelity solutions in large-scale systems. The basic principle of many ROMs is that the full order solution can be approximated using only a linear combination of 'dominant' modes. The procedure entails two stages: some typically expensive 'offline' calculations and inexpensive online calculations. In non-parametric transient cases, the Reduced Order Model (ROM) coefficients are only time dependent. In transient parametric cases, in the offline stage, the ROM is also trained on a set of different geometrical or physical parameter values, where as a result, a low-dimensional basis is computed. In the online stage, the ROM can be used for evaluating solutions on non-trained points in the parameter space or for time evolution if non-parametric PDEs are of interest. This process is achieved by projecting the full order equations (e.g. the Navier-Stokes equations) onto the calculated low-dimensional basis. An extensive review on projection-based ROMs can be found in \cite{Benner2015ASO}. 

In this work, the interest is focused on a hybrid POD - Galerkin - Proper Orthogonal Decomposition with Interpolation (PODI) based on Radial Basis Function method for parametric, time dependent PDEs governing turbulent fluid flow for non-isothermal mixing problems. POD-Galerkin has been applied to non-parametric PDEs with success across a wide range of applications. Beginning with the POD, this method was originally conceived as a method for compressing large datasets. An overview of POD can be found in \cite{cordier:hal-00417819}. In the field of fluid dynamics, Lumley \cite{lumley} was the first to apply POD to the study of turbulent flows. Later, Sirovich \cite{sirovich1987snap} proposed an 'optimised' method of performing the POD, known as the snapshot POD, where large data sets which are usually arising from large scale systems are available. Other applications of POD can be found in \cite{baltzer,bernero,rempfer}. 

The idea of PODI method, which was introduced by Buithanh in \cite{buithanh2003proper}, is based on the calculation of the POD coefficients using an interpolation technique. This method has been applied in aerodynamics in \cite{Dolci2016ProperOD} and in \cite{salmoiraghi} for parametrised geometries, as well as in cardiac mechanics in \cite{RAMA2016409}.     
  
The use of POD methods with Galerkin projection fall into the category of projection-based ROMs and has been applied in various areas for incompressible \cite{ravindran,Bourguet2011159,Lorenzi2016} and compressible flows \cite{Barone_stab_gal,Rowley2004115}.   
In regard to parametric PDEs, which are the main interest of this study, POD-Galerkin has been applied in \cite{ballarinmono} where the authors proposed a monolithic reduced order modelling approach for parametrised fluid-structure interaction problems. Also in \cite{Ballarin2015}, a stable POD-Galerkin method for the parametrised, incompressible steady Navier-Stokes equation is presented. Stabile \textit{et al} in \cite{Stabile2018} applied a POD-Galerkin ROM to the parametrised, incompressible unsteady Navier-Stokes equations. A POD-Galerkin model order reduction for parametric PDEs is also found in the study of haemodynamics in the work of Ballarin \textit{et al} \cite{BALLARIN2016609}. In \cite{saddamturbulent} and \cite{saddamdatadriven}, Hijazi \textit{et al} developed a POD-Galerkin ROM method for parametric Navier-Stokes equations in the turbulent regime while in \cite{StaBaZuRo2019} the authors focused on turbulence closure using VMS techniques. Another approach which uses Petrov-Galerkin projection for efficient nonlinear model order reduction is presented in \cite{Carlberg2010,xiao2013non,Carlberg2013623}.

In this work, the 3D, parametric, transient Navier-Stokes equations, weakly coupled with the energy equation are considered. In particular, the method presented in this work is applied to the modelling of turbulent thermal mixing in a T-junction pipe, a configuration commonly found in nuclear power reactor cooling systems and which plays a crucial role in the safety of reactors. The mixing of two different temperature streams leads to high transient temperature fluctuations in the pipe wall regions which can potentially lead to thermal fatigue and subsequent failure of the piping material (cracks formation, breakage etc). Turbulent thermal mixing has been studied both experimentally and computationally in \cite{FRANK20102313,WALKER2009116,NAIKNIMBALKAR20105901,AYHAN2012183,TUNSTALL2016672,Kuczaj2010}. In the computational case, various turbulent modelling techniques have been studied, including Large Eddy Simulation (LES) and Unsteady Reynolds Averaged Navier-Stokes (URANS). These methods, given the high Reynold's numbers and the nature of the problem, require fine 3D meshes which lead to very high CPU and memory costs. A ROM could therefore play a key role in such studies, giving the ability to obtain fast and realiable simulations. Model order reduction for nuclear applications has been previously applied in \cite{sartori2016reduced} for modelling the movement of the control rods in a nuclear reactor and in \cite{buchan2013pod} for reactor critically problems. This work extends the previous work of \cite{geostab} from laminar to turbulent heat transfer case. Other applications on POD-Galerkin ROM for coupled Navier-Stokes and energy equations include the work of Busto \textit{et al} \cite{BuStaRoCen2018} as well as in \cite{Raghupathy2009}. To the best of the authors knowledge, a ROM for modelling the parametric 3D Navier-Stokes equations, weakly coupled with the energy equation, including turbulence modelling is presented here for the first time. For the calculation of the reduced basis, two approaches are compared: a standard POD and a nested POD method. For the eddy viscosity term, an approach similar to the one developed in \cite{saddamturbulent,saddamdatadriven} is followed. This approach is based on the use of PODI with Radial Basis Function (RBF) interpolation for the calculation of the temporal coefficients of the eddy viscosity term. The POD-RBF method has been previously applied in the model order reduction context in \cite{dehghan2016proper} as a non-intrusive model order reduction method as well as in \cite{xiao2015non} for multiphase flow in porous media. In \cite{ostrowski}, a POD-RBF network for inverse heat conduction problems is presented.

This paper is organised as follows: in \ref{sec:math_form} the mathematical formulation is presented and in \ref{sec:ROM} the reduced order method is developed. In \ref{sec:training} a case of turbulent thermal mixing in T-junction is modeled and the ROM is trained and tested for several velocity inlet values. In \ref{sec:num_exp}, a comparison between standard and nested POD methods is performed on a particular case from \ref{sec:training}. In the final section, \ref{sec:concl} conclusions and perspectives are drawn, highlighting the directives for future improvements and developments. 
%%      MATH FORMULATION - sec:math_form
%%
%
%
%
%
% %%%%%%%%%%%%%%%%%%%%%%%%%%%%%%%%%%%%%%%%
 \section{Full Order Model - Mathematical Background}\label{sec:math_form}

In this section, the mathematical formulation of the FOM is presented. We consider the three dimensional, incompressible, transient, parametric Navier-Stokes equations weakly coupled with a three dimensional transient, parametric transport equation for heat. Since we are considering turbulent flows, some relevant background in turbulence modelling is presented. 

\subsection{Turbulence Modelling}
Turbulence is described by chaotic and random motions where the transported quantities (pressure, velocity etc) exhibit spatial and temporal fluctuations. Unlike laminar flows which can be directly numerically solved by many discretisation techniques (e.g. finite volume), practical modelling of turbulent flows involves additional approximations. While there is no exact delineation, internal flow (pipes or ducts) is usually considered turbulent if the non-dimensional Reynolds number, Re$>$4000, laminar if Re$<$2300 and transitional when 2300$<$Re$<$4000. The Reynolds number indicates the relative significance of the inertia forces in comparison to viscous forces. Turbulent motion greatly increases the effective diffusivity which leads to enhanced mixing, resulting in greater heat and momentum transfer.       

In the modelling of turbulent flows, the instantaneous velocity component $\bm{u}(\bm{x},t)$ is decomposed into a time-averaged component $\bm{\overline{u}}(\bm{x},t)$ and a fluctuating component $\bm{u'}(\bm{x},t)$, which can be expressed as:
\raggedbottom
\begin{equation}\label{eq:turbulence}
\bm{u}(\bm{x},t) = \bm{\overline{u}}(\bm{x},t) + \bm{u}'(\bm{x},t).
\end{equation}

The fluctuating component in the above equation (\ref{eq:turbulence}) is known as the turbulent fluctuation and is always three dimensional (in space), even for flows where the mean values change only in two dimensions. Turbulent flows contain rotational flow structures, the turbulent eddies, which vary in size. Larger eddies, acquire energy from the mean flow by a process called vortex stretching. The smaller eddies derive energy from the larger eddies through an energy cascade process and viscous dissipation converts turbulent energy from the smallest eddies (0.1-0.01mm - Kolmogorov microscales) into thermal internal energy leading to energy losses.

Regarding turbulence modelling, there are various computational approaches. These approaches include the Reynolds-Averaged Navier-Stokes method (RANS or URANS for the unsteady case),  Large Eddy Simulation (LES) and the most CPU expensive of all, Direct Numerical Simulation (DNS). RANS-URANS is the most widely used method in industry due to its relatively low computational cost in comparison with other methods for which the Navier-Stokes equations are solved in an averaged manner (ensemble or time). In LES the larger eddies are directly resolved whereas the smaller eddies (smaller than the mesh size) are modelled. LES has much greater computational cost than RANS and can be impractical for industrial applications. DNS numerically solves the Navier-Stokes equations and resolves the whole spectrum of eddies but this comes at vast computational cost.

In this work, URANS has been chosen due to its applicability in industry and reduced computational effort. Amongst the various methods available (Spalart-Almaras, $k-\epsilon$, $k-\omega$, Reynolds Stresses etc) in RANS/URANS, the $k-\omega$ SST is chosen due to its accuracy in the near wall region. There are two variations of the $k-\omega$ model, the standard $k-\omega$ and the Shear Stress Transport (SST) $k-\omega$. The former is described by a two-transport-equation model for $k$ and $\omega$ and the specific dissipation rate ($\epsilon/k$) is based on Wilcox \cite{wilcox1998turbulence}. The latter, is a combination of the standard $k-\omega$ model in proximity of the walls and the standard $k-\epsilon$ model in the bulk of the flow. To ensure smooth transition between the two different models, a blending function is used.                    

To derive the full-order equations which also include turbulence modelling, we start with the standard incompressible, transient parametric Navier-Stokes equations which are weakly coupled with the transient parametric heat transport equation. Considering an Eulerian framework and domain $\mathcal Q = \Omega \times [0,T_s] \subset \mathbb{R}^d\times\mathbb{R}^+$ with $d=1,2,3$, the equations are formulated as follow:  
%%    STRONG FORM NAVIER STOKES - eq:navstokes
%%
\begin{equation}\label{eq:navstokes}
\begin{cases}
\frac{\partial \bm{u}(\bm{x},\bm\mu,t)}{\partial t}+ \bm{\nabla} \cdot (\bm{u}(\bm{x},\bm\mu,t) \otimes \bm{u}(\bm{x},\bm\mu,t))- \\
-\bm{\nabla} \cdot \nu\bm{\nabla} \bm{u}(\bm{x},\bm\mu,t)=-\bm{\nabla}p(\bm{x},\bm\mu,t) &\mbox{ in } \mathcal Q,\\
\bm{\nabla} \cdot \bm{u}(\bm{x},\bm\mu,t)=\bm{0} &\mbox{ in } \mathcal Q,\\
\frac{\partial \theta(\bm{x},\bm\mu,t)}{\partial t} + \nabla \cdot (\bm{u}(\bm{x},\bm\mu,t) \theta(\bm{x},\bm\mu,t) ) -\\
- \alpha_{dif}\Delta \theta(\bm{x},\bm\mu,t) = 0 &\mbox{ in } \mathcal Q,\\
\bm{u} (\bm{x},\bm\mu) = \bm{f}(\bm{x},\bm\mu) &\mbox{ on } B_{In},\\
\theta (\bm{x},\bm\mu) = g(\bm{x}) &\mbox{ on } B_{In},\\ 
\nabla\theta(\bm{x},\bm\mu)\cdot\bm{n}= 0 &\mbox{ on } B_{w},\\
\bm{u} (\bm{x},\bm\mu) = \bm{0} &\mbox{ on } B_{w},\\ 
\nu(\nabla \bm{u} - p\bm{I})\bm{n} = \bm{0} &\mbox{ on } B_{o},\\ 
\bm{u}(0,\bm{x})=\bm{h}(\bm{x}) &\mbox{ in } T_0,\\     
\theta(0,\bm{x})=e(\bm{x}) &\mbox{ in } T_0,\\        
\end{cases}
\end{equation}  
where $\bm{u}$ is the fluid velocity, $p$ the normalised pressure, $\theta$ is the fluid temperature, $\alpha_{dif}$ is the thermal diffusivity and $\nu$ is the kinematic viscosity. The vector of parameters is represented by the greek letter $\bm\mu$. $B_{In} = \Gamma_{In} \times [0,T_s]$, $B_w = \Gamma_{w} \times [0,T_s]$ and $B_o = \Gamma_o \times [0,T_s]$, $T_s$ represents the time-period of the simulation, $\Gamma = \Gamma_{In} \cup \Gamma_{w} \cup \Gamma_{o}$ is the boundary of $\mathcal{Q}$ and consists of three different parts $\Gamma_{In}$, $\Gamma_{o}$ and $\Gamma_w$ that indicate, respectively, inlet(s), outlet and physical wall boundaries (\ref{fig:mesh}. The functions $\bm{f}(\bm{x},\bm\mu)$ and $g(\bm{x})$ represent the boundary conditions for the non-homogeneous boundaries. $\bm{h}(\bm{x})$ and $e(\bm{x})$ denote the initial conditions for velocity and temperature at $t=0$. Time independence of the boundary conditions $\bm{f}$ and $g$ is also assumed. Since we are dealing with the incompressible equations, the density $\rho$ has been omitted and is assumed to have the value of $1$. In this work, the parametric dependency is on the velocity inlet boundary conditions.

Similar to equation (\ref{eq:turbulence}) the instantaneous velocity, pressure and temperature can be written as:
\raggedbottom
\begin{eqnarray}\label{eq:turbulenceall}
\bm{u}(\bm{x},t) = \bm{\overline{u}}(\bm{x},t) + \bm{u}'(\bm{x},t)\\
p(\bm{x},t) = \overline{p}(\bm{x},t) + p'(\bm{x},t),\\
\theta(\bm{x},t) = \overline{\theta}(\bm{x},t) + \theta'(\bm{x},t) .
\end{eqnarray}

When using RANS modelling, the mean components of the fluctuating terms are zero. Taking this into account and substituting equations (\ref{eq:turbulenceall}) into equations (\ref{eq:navstokes}), the Reynolds-averaged Navier-Stokes and energy equations are:
\raggedbottom
\begin{equation}
\begin{cases}
\nabla\cdot\bm{\overline{u}} = 0, \\
\frac{\partial\bm{\overline{u}}}{\partial t} + \nabla\cdot(\bm{\overline{u}}\otimes\bm{\overline{u}}) - \nu\Delta\bm{\overline{u}} = -\nabla\overline{p} - \nabla\cdot\bm{\tau}^R, \\
\frac{\partial\overline{\theta}}{\partial t} + \nabla\cdot(\bm{\overline{u}}\overline{\theta}) - (\alpha_{{dif}}\Delta\overline{\theta}) - \nabla\cdot\bm{H}^R = 0 ,
\end{cases}
\end{equation}

where the extra terms $\bm{\tau}^R = -(\overline{\bm{u}'\bm{u}'})$ and $\bm{H}^R = (\overline{\bm{u}\theta} - \overline{\bm{u}}\overline\theta)$ are the Reynolds stress tensor and the heat flux term respectively. Therefore, these two new terms must be modelled in order to close the system of the equations. One possible solution, which we will follow here is the eddy viscosity models with the Boussinesq hypothesis. The Reynolds stress is modelled using an eddy viscosity $\nu_t$ 

\raggedbottom
\begin{equation}
\bm{\tau}^R = -(\bm{\overline{u'u'}}) = 2\nu_t\bm{\overline{D}}^R - 2/3k\bm{I} = \nu_t[\nabla\bm{\overline{u}}+(\nabla\bm{\overline{u}})^T]
- 2/3k\bm{I},  
\end{equation}
 
while the heat flux using a gradient diffusion hypothesis as follows:

\raggedbottom
\begin{equation}
\bm{H} = -\alpha_{{dif}_t}\nabla\overline{\theta},
\end{equation}

where $\alpha_{{dif}_t}$ is the turbulent thermal diffusivity defined as $\nu_t/Pr_t$ and $Pr_t$ the turbulent Prandtl, $\bm{\overline{D}}^R$ the Reynolds-averaged strain rate tensor and $\bm{I}$ the identity matrix.
 
% where $\alpha_{{dif}_{eff}} = \nu/Pr + \nu_t/Pr_t$.   

Combining all of the above and denoting from now on $\bm{\overline{u}}$, $\overline{p}$ and $\overline\theta$ as $\bm{u}$, $p$ and $\theta$ respectively, equations (\ref{eq:navstokes}) including also the turbulence modelling are as follows: 

\begin{equation}\label{eq:navstokesdiff}
\begin{cases}
\frac{\partial \bm{u}}{\partial t} + \nabla\cdot(\bm{u}\otimes\bm{u}) = \nabla \cdot \big[(\nu + \nu_t)\cdot \frac{1}{2}\big(\nabla\bm{u} + \\
+(\nabla\bm{u}^T)\big) -p\bm{I}\big], \\
\bm{\nabla} \cdot \bm{u}=\bm{0}, \\
\frac{\partial \theta}{\partial t} + \nabla \cdot (\bm{u} \theta ) - (\alpha_{dif}+\alpha_{dif_t})\Delta \theta = 0, \\
\nu_t = j(k,\omega), \\
\bm{u}  = \bm{f} &\mbox{ on } B_{In},\\
\theta  = g &\mbox{ on } B_{In},\\
\bm{u} = \bm{0} &\mbox{ on } B_{w},\\ 
\nabla\theta\cdot\bm{n}= 0 &\mbox{ on } B_{w},\\ 
\nu(\nabla \bm{u} - p\bm{I})\bm{n} = \bm{0} &\mbox{ on } B_{o},\\ 
\bm{u}(0)=\bm{h} &\mbox{ in } T_0,\\     
\theta(0)=e &\mbox{ in } T_0,\\  
\end{cases}
\end{equation}  

within the space-time domain $\mathcal Q$. The turbulent kinetic energy is $k$, $k=\frac{1}{2}\bm{u}\cdot\bm{u}$ and $\omega = \frac{\epsilon}{k\beta^{*}}$ is the specific dissipation rate with $\beta^{*}$ being a constant of proportionality and $\epsilon \propto \frac{\partial k}{\partial t}$ is the dissipation rate. For the $k$ and $\omega$ transport equations the reader should refer to \cite{Versteeg1995AnIT}. For brevity, the parametric, spatial and temporal dependency on the momentum, continuity and energy equations has been omitted.

\subsection{Finite Volume Approximation}
\begin{figure*}
\centering
\includegraphics[width=0.26\textwidth]{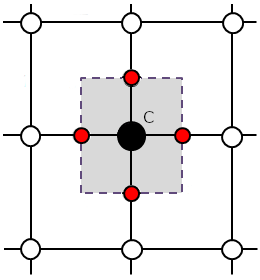}
\caption{A 2D control volume with C being the centroid of the volume.}\label{fig:volume}
\end{figure*}

The full order partial differential equations (\ref{eq:navstokesdiff}) are transformed into a discrete system of algebraic equations using the finite volume method. The solution domain is discretised with the generation of a computational mesh in a number of arbitrary cells or control volumes over which the equations are solved (spatial discretisation) and the variables of interest are calculated at the centroid of the control volumes (\ref{fig:volume}). The time is also discretised in a number of timesteps (temporal discretisation). The integral form of the equations are discretised over each cell so that, at the discrete level, the mass and momentum are conserved.     

To give a brief description of the discretisation process, a general transport equation for a transported quantity $\phi$, in a control volume $\Omega_i$, where the density $\rho$ is taken equal to $1$ and no sources are considered, can be expressed as:
\raggedbottom
\begin{equation}\label{eq:transport}
\int_{\Omega_i} \frac{\partial\varphi}{\partial t}d\Omega + \int_{\Omega_i} \nabla\cdot(\bm{u}\varphi)d\Omega - \int_{\Omega_i} \nabla\cdot(\Gamma_\varphi\nabla\varphi)d\Omega = 0 ,
\end{equation}

where $\Gamma_\varphi$ is the diffusion coefficient of the transported quantity.  The above transport equation (\ref{eq:transport}) consists of 3 terms with the first one being the temporal term, the second one the convective term and the last one the diffusive term. As this equation is second order, a second or higher order discretisation scheme is needed. The volume integrals in the convective and diffusive terms are converted into surface integrals using Gauss theorem, which for a vector $\bm{w}$ is:
\raggedbottom
\begin{equation}
\int_{\Omega_i} \nabla\cdot \bm{w}d\Omega = \oint_{\partial\Omega_i} \bm{n}dS\cdot\bm{w} ,
\end{equation} 

where $\bm{n}$ denotes the normal vector pointing outward of the bounding surface $\partial\Omega_i$. The volume and surface integrals are then approximated by a discrete sum over the faces of the control volume $\Omega_i$. The values on the faces of the control volume are calculated by interpolating the values from the adjacent control volumes. The most widely used interpolation schemes are the central differencing and upwind differencing schemes. Regarding the temporal term, there is also a variety of differencing schemes including the backward differencing scheme, Crank-Nicolson, Euler etc. For a detailed review on the finite volume discretisation the reader could refer to \cite{Moukalled:2015:FVM:2876154}.

\section{Reduced Order Model Formulation}\label{sec:ROM}

A POD-Galerkin method is used for the development of the ROM. We adopt a procedure similar to that used for the laminar case in \cite{geostab}, combined with the procedure first introduced in \cite{saddamturbulent} for the modelling of the eddy viscosity term when high Reynolds numbers are considered. This method is extended to the modelling of the eddy diffusivity when the temperature transport equation is considered. 

One can decompose the velocity, pressure and temperature fields as a linear combination of basis functions multiplied by temporal coefficients as follows:
\raggedbottom
\begin{align}
\bm{u}(\bm{x},\bm{\mu},t) \approx\bm{u}_s = \sum_{i=1}^{N_{u}^s}\alpha_i(\bm{\mu},t)\bm{\phi}_i(\bm{x}), \\
p(\bm{x},\bm{\mu}, t) \approx p_s = \sum_{i=1}^{N_p^s}b_i(\bm{\mu},t)\psi_i(\bm{x}), \\
\theta(\bm{x},\bm{\mu}, t) \approx \theta_s =  \sum_{i=1}^{N_{\theta}^s}c_i(\bm{\mu},t)\chi_i(\bm{x}) ,
\end{align}

where $\bm{\phi}(\bm{x})$, $\psi(\bm{x})$ and $\chi(\bm{x})$ are spatial basis functions and $\alpha(\bm{\mu},t)$, $b_i(\bm{\mu},t)$ and $c_i(\bm{\mu},t)$ are the temporal coefficients for the velocity, pressure and temperature fields respectively. $N^s_{u}$, $N^s_{p}$ and $N^s_{\theta}$ is the cardinality of the POD spaces for velocity, pressure and temperature respectively. Considering the FOM, which exhibits parametric and temporal dependency, is solved for $N_p$ training points of a parameter of interest $\mu^p \in \mathcal P = \{\mu^1,\dots ,\mu^N_p\}$ and for $N_t$ time instances, $t^k\in\{t^1,\dots, t^k\}\subset [0,T_s]$. Therefore, the total number of velocity snapshots is $N_u^s = N_p\times N_t$.  

For turbulent modelling, a new term is introduced which approximates the eddy viscosity term:
\raggedbottom
\begin{equation}\label{eq:eddyvisc}
\nu_t(\bm{x},\bm{\mu},t) \approx\nu_{t_s}=\sum_{i=1}^{N_{\nu_t}^s}l_i(\bm{\mu},t)\xi_i(\bm{x}),
\end{equation}
  
where $\xi(\bm{x})$ and $l(\bm{\mu},t)$ represent the basis functions and the temporal coefficients for the eddy viscosity field. 

The basis functions are calculated using the POD method for parametric problems, where an equispaced grid sampling over a range of parameters and time is performed. This results to a global snapshot matrix which contains both parameter and time information. Therefore, for the velocity, the snapshot matrix, $\bm{U}_s$, is formed as
\raggedbottom
\[
\bm{U}_s = \begin{bmatrix}\label{eq:globalmatrix} 
    u_{1}^{1}(\mu^1) & \dots & u_{1}^{N_{t}}(\mu^1) &  \dots & u_{1}^{N_{t}}(\mu^{N_p}) \\
    \vdots & \dots & \vdots & \dots & \vdots & \\
    u_{N_{u}^h}^{1}(\mu^1) & \dots & u_{N_{u}^h}^{N_{t}}(\mu^1) & \dots & u_{N_{u}^h}^{N_{t}}(\mu^{N_p}) 
    \end{bmatrix},
\]

where $N_{u}^h$ are the number of degrees of freedom which is equal to the number of grid points times the number of the components. 
For any vector or scalar function $\bm{w}$, the POD method seeks to find an orthonormal basis $\bm\phi(\bm{x})$ of rank $N_w^r$ which minimises the quantity:

\begin{small}
\begin{equation}\label{eq:pod_energy}
\mbox{arg min}\sum_{i=1}^{N_{w}^s}||\bm{w}_{i} - \sum_{j=1}^{N_{w}^r}\ltwonorm{\bm{w}_{i}(\bm{x},\bm{\mu},t),\bm{\phi}_j(\bm{x})}\bm{\phi}_j(\bm{x}) ||^2_{L^2},
\end{equation}
\end{small}

with,

\begin{equation}
\ltwonorm{\phi_{i},\phi_{j}} = {\delta}_{{ij}}, \forall i,j = 1,\dots,N^{s}_{w} .
\end{equation}

The vector function $\bm{w}$ can be the velocity vector $\bm{u}(\bm{x},\bm\mu,t)$ or the scalar functions for temperature $\theta(\bm{x},\bm\mu,t)$ and pressure $p(\bm{x},\bm\mu,t)$. The above problem (\ref{eq:pod_energy}) is equivalent to solving the eigenvalue problem:
\raggedbottom
\begin{eqnarray}\label{eq:pod_eigenvalue}
\bm{C}\bm{W} &=& \bm{W}\bm\lambda , \\
{C}_{ij} &=& \ltwonorm{\bm{u}_i,\bm{u}_j}, \forall i,j = 1,\dots,{N}^{s}_{u} ,  
\end{eqnarray} 

where $\bm{C}$ is the correlation matrix, $\bm{W}$ a matrix with the eigenvectors and $\bm\lambda$ a diagonal matrix with the eigenvalues. 

The original POD basis, $\mathcal{V}_{POD}=$span$[\bm{\phi}_1,\bm{\phi}_2,...\bm{\phi}_{N_u^s}]$ can be truncated according to the following energy criterion:
\raggedbottom
\begin{equation}\label{eq:energycr}
E_{N^r_u}=\frac{\sum_{i=1}^{N^r_u}\lambda_{i}}{\sum_{j=1}^{N_{u}^s}\lambda_{j}} ,
\end{equation}

where $\lambda_{i}$ are the eigenvalues and $N^r_u << N^s_{u}$ denotes the cardinality of the reduced POD space $\mathcal{V}_{POD_r}$.

The orthogonal POD basis functions, $\bm{\phi}(\bm{x})$, are calculated and normalised as:
\raggedbottom
\begin{eqnarray}\label{eqn2}
\bm{\phi}_j &=&   \frac{1}{\sqrt{\lambda}_{i} N^r_u}\sum_{j=1}^{N^r_u}\bm{u}_j \bm{W}_{ij}, \\
\ltwonorm{\bm{\phi}_i,\bm{\phi}_j} &=& \delta_{ij}  \hspace{0.5cm}\forall\mbox{ } i,j = 1,\dots,N^r_u.
\end{eqnarray}

Other methods for the generation of the snapshot matrix include, for example, Greedy POD \cite{Haasdonk2013,urban2014greedy}, the goal-oriented POD-greedy sampling \cite{hoang2015efficient}, the Proper Generalised Decomposition (PGD) \cite{ChinestaEnc2017}, the nested POD \cite{doi:10.1002/gamm.201610011} as well as other approaches, \cite{Rozza2008,ChinestaEnc2017,Kalashnikova_ROMcomprohtua,quarteroniRB2016,Chinesta2011,Dumon20111387}. In this work, the nested POD will tested against the standard POD method. 

\subsection{Nested POD}\label{subsec:nested}
The standard POD method becomes too expensive for training spaces with many training points as this leads to large dense matrix eigenproblems. Considering for example the velocity snapshots, the computational effort need for the solution of the quadratic eigenvalue problem, (\ref{eq:pod_eigenvalue}), scales with $\mathscr{O}([N_{u}^s]^3)$. The nested POD method approximates the global POD space by solving one small eigenvalue problem for each local parameter space. The local POD bases are then weighted by the eigenvalues and a global snapshot matrix is created by appending the local weighted POD basis. A standard POD is then performed on the global snapshot matrix and the basis is calculated according to the method described in (\ref{sec:ROM}). The advantage of the nested POD is the numerical efficiency over the standard POD method since the computational effort for the former scales as $\mathscr{O}([N_{t}^3\cdot N_p +[N_{u}^{\emph{nested}}]^3)$, where $N_{u}^{\emph{nested}}$ is the dimension of the global snapshot matrix resulting from the nested POD.

The nested POD acts on the following $N_p$ local matrices:

\[
\bm{U}_{\emph{nested}}^i = \begin{bmatrix}\label{nestedmatrix}
    u_{1}^{1}(\mu^i) & u_{1}^{2}(\mu^i) & \dots & u_{1}^{N_{t}}(\mu^i) \\
    \vdots & \vdots & \dots & \vdots & \\
    u_{N_{u}^h}^{1}(\mu^i) & u_{N_{u}^h}^{2}(\mu^i) & \dots & u_{N_{u}^h}^{N_{t}}(\mu^i)
    \end{bmatrix} ,
\]

where $i$ runs from $1$ to $N_p$. A POD is then performed on each local matrix $\bm{U}_{\emph{nested}}^i$ and we obtain corresponding local bases.

The global POD matrix is then formed by selecting the first $N_{t}^{rn}$ modes from each local POD basis according to the energy criterion (\ref{eq:energycr}), weighted by their eigenvalues. The concatenation of the weighted modes results in the following global POD matrix: 
\raggedbottom
\[
\bm{U}_{s_n} = \left[\begin{smallmatrix}\label{smallglobal}
    \lambda_{1}^{1}\phi_{1}^{1}(\mu^1) & \dots & \lambda_{1}^{N_{t}^{rn}}\phi_{1}^{N_{t}^{rn}}(\mu^1) & \dots & \lambda_{1}^{1}\phi_{1}^{N_{t}^{rn}}(\mu^{N_p}) \\
    \vdots & \dots & \vdots & \dots & \vdots & \\
    \lambda_{1}^{1}\phi_{N_u^h}^{1}(\mu^1) & \dots & \lambda_{1}^{N_{t}^{rn}}\phi_{N_{u}^h}^{2}(\mu^1) & \dots & \lambda_{1}^{N_{t}^{rn}}\phi_{N_{u}^h}^{N_{t}^{rn}}(\mu^{N_p})
    \end{smallmatrix}\right].
\]
 
Once the snapshot matrix $\bm{U}_{s_n}$ is formed, the procedure is the same as that described in section (\ref{sec:ROM}).

\subsection{Galerkin Projection}   
\raggedbottom
For the calculation of the temporal coefficients, a Galerkin approach is followed for the velocity, pressure and temperature fields, where the original equations are projected onto the reduced basis. The momentum equation is projected onto the spatial POD basis $\bm{\phi}(\bm{x})$ while the continuity equation is projected onto the pressure spatial basis $\psi(\bm{x})$, using the supremizer method \cite{Rozza2007,Ballarin2014}. The energy equation is projected onto the temperature spatial basis $\chi(\bm{x})$. The projection results in a set of ordinary differential equations showing the evolution of the temporal coefficients for velocity, pressure and temperature:
\raggedbottom
\begin{align}\label{eq:reducedeq}
\bm{M}\bm{\dot\bm{\alpha}} - \bm\alpha ^T \bm{Q}\bm{\alpha} + (\nu+\nu_t)\bm{l}^T\big(\bm{Q}_{T1} + \bm{Q}_{T2}\big)\bm{\alpha} - \bm{P}\bm{b} &=& 0, \\
\bm{R}\bm{\alpha} &=& 0 , \\
\bm{K}\bm{\dot{c}} - \bm{\alpha} ^T \bm{G}\bm{c} - \alpha_{dif_t}\bm{N}\bm{c} &=& 0 ,
\end{align}

where the dotted terms $\dot{\bm{\alpha}}$ and $\dot{\bm{c}}$ represent time derivatives and the reduced matrices are:

\raggedbottom
\begin{align}
(\bm{M})_{{ij}} & =  \ltwonorm{\bm{\phi}_i, \bm{\phi}_j}, \\
(\bm{Q})_{{ijk}} & =  \ltwonorm{\nabla\cdot(\bm{\phi}_i\otimes\bm{\phi}_j),\bm{\phi}_k}, \\
(\bm{Q})_{T1_{ijk}} & \ltwonorm{\xi_{i}\Delta\bm{\phi}_j,\bm{\phi}_k}, \\
(\bm{Q})_{T2_{ijk}} & \ltwonorm{\nabla\cdot\xi_{i}(\nabla\bm{\phi}^T_j),\bm{\phi}_k}, \\
(\bm{P})_{{ij}} & =   \ltwonorm{\nabla\psi_i,\bm{\phi}_j},\\
(\bm{R})_{{ij}} & =   \ltwonorm{\nabla\cdot\bm\phi_i,{\psi}_j}, \\
(\bm{K})_{{ij}}  & =  \ltwonorm{\chi_i, \chi_j}, \\
(\bm{G})_{{ijk}} & =  \ltwonorm{\nabla\cdot(\bm{\phi}_i\chi_j),\chi_k}, \\
(\bm{N})_{{ij}}  & =  \ltwonorm{\Delta\chi_i,\chi_j} .
\end{align}

\subsection{Radial Basis Functions}
 
Considering the eddy viscosity field, the vector of temporal coefficients $\bm{l}(t)$ is computed using the non-intrusive method of Radial Basis Function (RBF) interpolation according to the procedure described in \cite{saddamturbulent}. In this way, since there is no projection of the eddy viscosity modes onto the turbulence modelling equations, the ROM is independent of the turbulent method which are used in the FOM.     

The temporal coefficients, $l(\bm\mu,t)$, for a new value of the parameter, are calculated during the online stage as a linear combination of $N_{\nu_t}^s$ ($j=1,2,\dots,N_{\nu_t}^s$) chosen radial basis functions kernels $\bm{\Theta}$:
\raggedbottom
\begin{equation}\label{eq:rbf} 
l_{i}(\bm\mu,t) = \sum_{j=1}^{N_{\nu_t}^s}w_{i,j}\Theta_{i,j}(||\bm\mu-\bm\mu_{j} ||_{L2}),
\mbox{ for } i = 1,2,\dots,{N}^{r}_{\nu_t},  
\end{equation} 

where, $\bm{w}$ represents the vector of the linear weights, $\bm\mu$ are the sampling points (centers) corresponding to eddy viscosity snapshots $\nu_{t}$ and $\bm\mu$ is the value of the new input parameter which does not coincide with any of the training points. For the RBFs, $\bm\Theta$, various kernels can be used. In this work, Gaussian kernels are considered, defined as follow:
\raggedbottom
\begin{equation}
\bm\Theta(||\bm\mu-\bm\mu_{j} ||_{L2}) = exp(-\gamma||\bm\mu-\bm\mu_{j} ||^2_{L2}) ,
\end{equation}

where $\gamma$ is a parameter which determines the spread of the kernel. The RBF monotonically decreases as we move away from the centre. Gaussian RBFs response is local, meaning that their response is the best in the area near to the centre, in contrast to multiquadratic RBFs which are global. The unknowns are computed in the offline stage. 

The vector of the weights is calculated by solving the following linear system:
\raggedbottom
\begin{equation}\label{eq:weights}
\sum_{j=1}^{N_{\nu_t}^s}w_{i,j}\Theta_{i,j}(||\bm\mu_{k}-\bm\mu_{j} ||_{L2}) =  l(\bm\mu,t)_{i,k}, \\
\mbox{ for } k = 1,2,\dots,{N}^{s}_{\nu_t}.
\end{equation} 
 
The above equation can be written in matrix form and solved for the unknown weights as:

\begin{equation}
\bm\Theta\bm{w}=\bm{l} \Leftrightarrow \bm{w} = \bm\Theta^{-1}\bm{l},
\end{equation}

provided that the matrix $\bm\Theta$ is non-singular, therefore invertible and the temporal coefficients, $l(\bm\mu,t)_{i,j}$, are calculated by projecting the eddy viscosity snapshots onto the spatial eddy viscosity modes as:
\raggedbottom
\begin{equation}\label{eq:eddymodes}
l_{i,j}(\bm\mu,t) = \ltwonorm{\xi_{i},\nu_{t_j}}, \\
\mbox{ for } j = 1,2,\dots,{N}^{s}_{\nu_t} .
\end{equation} 
 
Once the unknown weights are calculated, the system of equations (\ref{eq:reducedeq}), is ready to be solved.

\subsection{Treatment of Boundary Conditions}

The interest in the present work is the construction of a ROM with parametric boundary conditions. This means that, given a proper sampling of a given parameter, the ROM would be capable of simulating cases for new, untrained values of the chosen parameter provided these values belong to the training manifold. To achieve this property, a similar procedure to that introduced in \cite{NME:NME537} is followed. 

This method first employs a control function to make the snapshots homogeneous and independent of the boundary conditions by subtracting the inhomogeneous boundary value and then, at the reduced order level, the new boundary value is enforced and the contribution of the control function values is added back. Taking as example the velocity snapshots, these are homogenised as :
\raggedbottom
\begin{equation}
\bm{u}'(\bm{x},\bm\mu,t) = \bm{u}(\bm{x},\bm\mu,t) - \sum_{j=1}^{N_{BC}}u_{D_j}(\bm\mu,t)\bm{\zeta}_{c_j} ,
\end{equation} 

where $N_{BC}$ is the number of parametrised boundary conditions and $u_D$ are scaling coefficients that make the snapshots homogeneous. The control functions, $\bm{\zeta}$, are calculated by assigning:

\begin{equation}
    \bm{ \zeta_c} = 
    \begin{cases}
        1, & \text{at the Dirichlet boundary } \\
        0, & \text{on the other boundaries} ,
    \end{cases}
\end{equation}
      
and then are computed as $\bm{\zeta}_{c} = \frac{1}{N_{u}^s} \sum_{l=1}^{N_{u}^s} \bm{u}_{i}$ by solving a full order model. The POD is then applied to the homogeneous snaphots and the velocity field is approximated as:
\raggedbottom
\begin{equation}
\bm{u}(\bm{x},\bm\mu,t) = \sum_{j=1}^{N_{BC}}u_{D_j}(\bm\mu,t)\bm{\zeta}_{c_j} + \sum_{i=1}^{N_{u}^s}\alpha_i(t,\bm\mu)\bm{\phi}_i(\bm{x}).
\end{equation}

A similar procedure is followed for the temperature parametrised boundary conditions. For a detailed discussion the reader should refer to \cite{Stabile2017,geostab}.

\section{Applications}

The mathematical framework described in the previous sections is tested on a 3D T-junction pipe, where different temperature water streams are mixing in the tee area at high Reynolds numbers. The boundary of the domain $\Omega$, denoted with $\Gamma$, consists of four parts $\Gamma = \Gamma_{m} \cup \Gamma_{b} \cup \Gamma_{w} \cup \Gamma_{o}$ as shown in figure \ref{fig:mesh}. The geometrical properties of the pipe are shown in table \ref{tab:parameters} while the total length of the pipe is $L=3.0$m. In \ref{sec:training} the ROM is trained on ten different sets of inlet velocity values (table \ref{tab:sampling_points}) and then is tested on four sets of values within the range of the training space (table \ref{tab:test_points}). Then, in \ref{sec:num_exp}, the set D is selected and a comparative study between standard and nested POD methods is performed.      

\begin{table*}[!tbp]  % title of Table
\centering % used for centering table
\begin{tabular}{ l | c | c }
\hline\hline
& Main Pipe & Branch Pipe  \\ [0.5ex]
\hline
$u$ (m/s) & values in table (\ref{tab:test_points}) & values in table (\ref{tab:test_points})  \\ \hline
$\theta$ (K) & 292.15 & 309.5 \\ \hline
$D$ (m) & 0.14 & 0.08 \\ \hline
$Re$  & values in table (\ref{tab:test_points}) & values in table (\ref{tab:test_points}) \\
\hline
\end{tabular}
\caption{Summary of the physical parameters for the reduced order model.}\label{tab:parameters}
\end{table*}

\subsection{Numerical study: Thermal mixing in T-junction pipe}\label{sec:training}

A parametric turbulent case is studied, where the velocity inlet boundary conditions on both inlets is parametrised. This example has been chosen as the parametric response is non-linear and a proper training is needed. The training space is constructed using ten sets of sampling points for the two velocity inlets, as shown on table (\ref{tab:sampling_points}). The reduced order model is evaluated on four sets of values (table \ref{tab:test_points}) which belong to the training range but they do not coincide with the samples used to train the ROM. It is known from the literature \cite{FRANK20102313,FENG2018275,TUNSTALL2016170} that the modeling of turbulent non isothermal mixing in T-junctions appears to be a challenging testcase and the Large Eddy Simulation (LES) method is usually required for the simulation of transient velocity and temperature fields. However, since the intention of this work is to compare two different numerical approaches, the full order (FOM) simulation performed in OpenFOAM \cite{OF} and the reduced order (ROM) simulation performed in ITHACA-FV C++ library \cite{RoSta17}, the  URANS $k$-$\omega$ SST turbulence model has been used for time-saving reasons. As the hybrid PODI-Galerkin method presented in this paper is non-intrusive with respect to the eddy viscosity field, any type of turbulence modelling could be in principle used for the construction of the ROM.

\begin{figure*}[!tbp]  
  \centering 
  \begin{minipage}[b]{0.45\textwidth}
    \includegraphics[width=\textwidth]{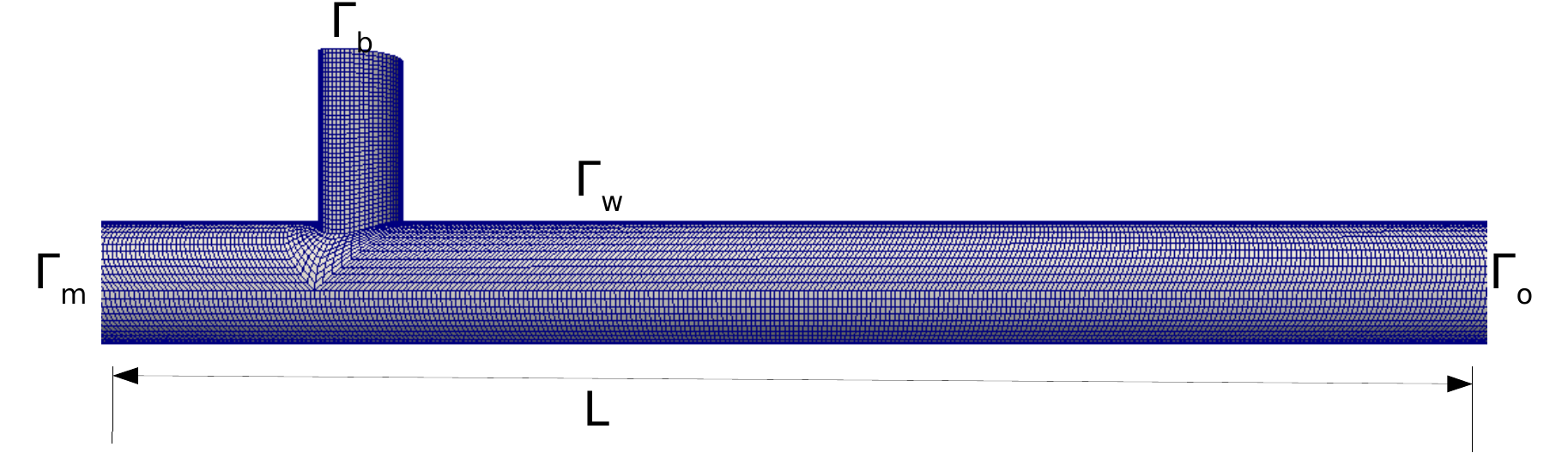}
    \caption{Computational mesh of the T-junction pipe.}\label{fig:mesh}
  \end{minipage}
  \hfill
  \begin{minipage}[b]{0.40\textwidth}
    \includegraphics[width=\textwidth]{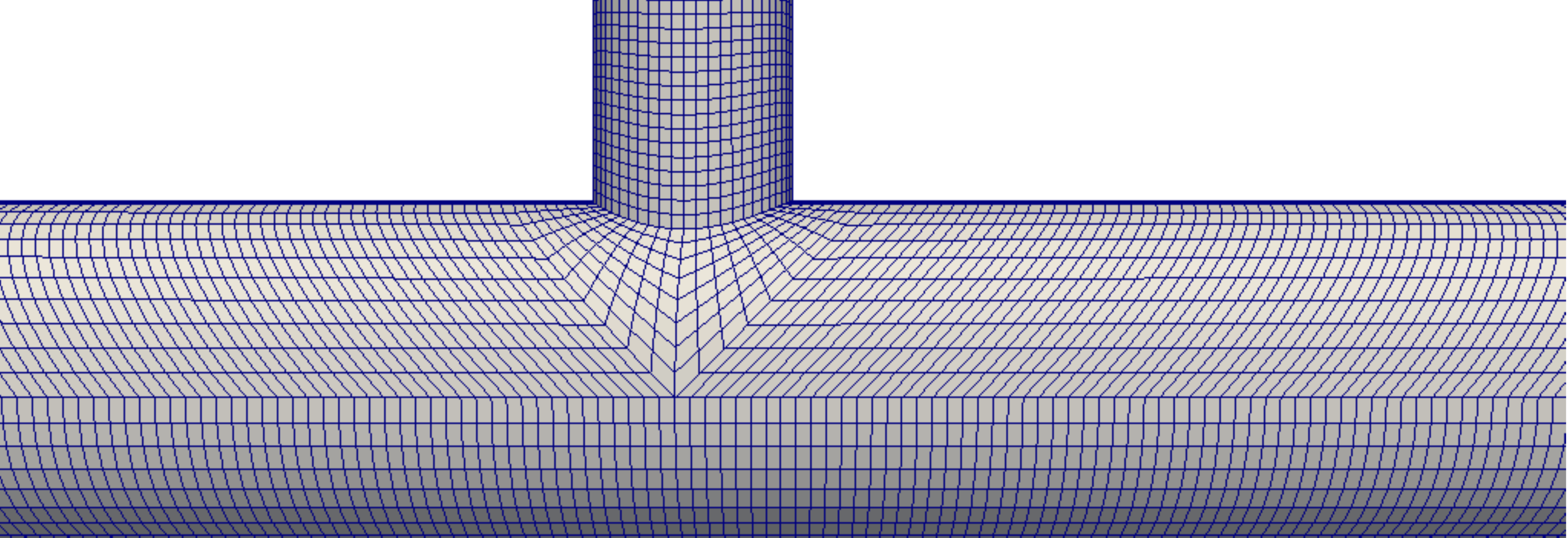}
    \caption{Mesh in T-junction region.}\label{fig:mesh_t}
  \end{minipage} 
\end{figure*}

In the offline phase, the FOM is modelled using the transient pisoFoam solver \cite{ISSA198640}. The simulation time is set to $3$s with timestep $dt = 0.0025$s so that the flow reaches the outlet. The computational mesh, shown in figures (\ref{fig:mesh}) and (\ref{fig:mesh_t}), consists of 291816 cells (hexahedral). The $y^+$ value of the mesh is 126 and for this reason wall functions have been used (kqRWallFunction for $k$ and omegaWallFunction for $\omega$ \cite{chalmers}). The spatial and temporal discretisation schemes are summarised in table (\ref{tab:schemes}). The FOM is run 10 times (one for each sampling pair) with snapshots being collected every $0.1$s. This makes a total of $30\cdot 10=300$ snapshots ($30$ snapshots per sampling pair) for each field. Therefore, $N_u^s=N_p^s=N_{sup}^{s}=N_{\theta}^s=N_{\nu_t}^s=300$. The reduced basis is computed with the POD method. The POD is applied directly on the global snapshots matrices which contains both parameter and time information in an equispaced grid, according to (\ref{eq:globalmatrix}) and the reduced spaces are then chosen according to the decay of the eigenvalues shown in figure (\ref{fig:energyeigen}). The retained modes are $6$ for velocity, $10$ for pressure and supremizer, $10$ for the eddy viscosity and $11$ modes for temperature.

\begin{figure*}[!tbp]
  \centering 
  \begin{minipage}[b]{0.45\textwidth}
    \includegraphics[width=\textwidth]{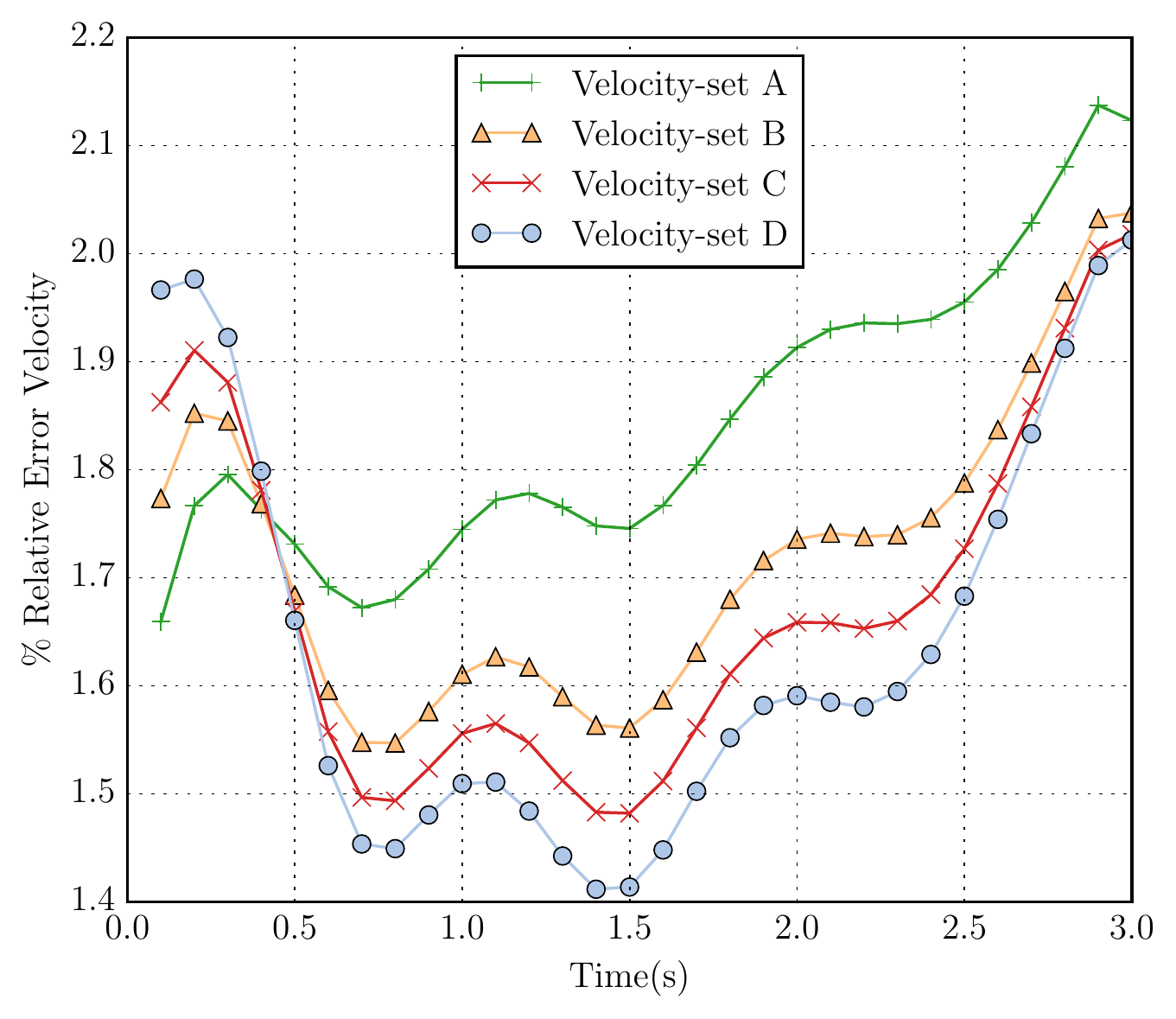}
    \caption{$\% \epsilon_{L^2}(t)$ error of velocity field for the four test sets (\ref{tab:test_points}).}\label{fig:errorl2vel}
  \end{minipage}
  \hfill
  \begin{minipage}[b]{0.45\textwidth}  
    \includegraphics[width=\textwidth]{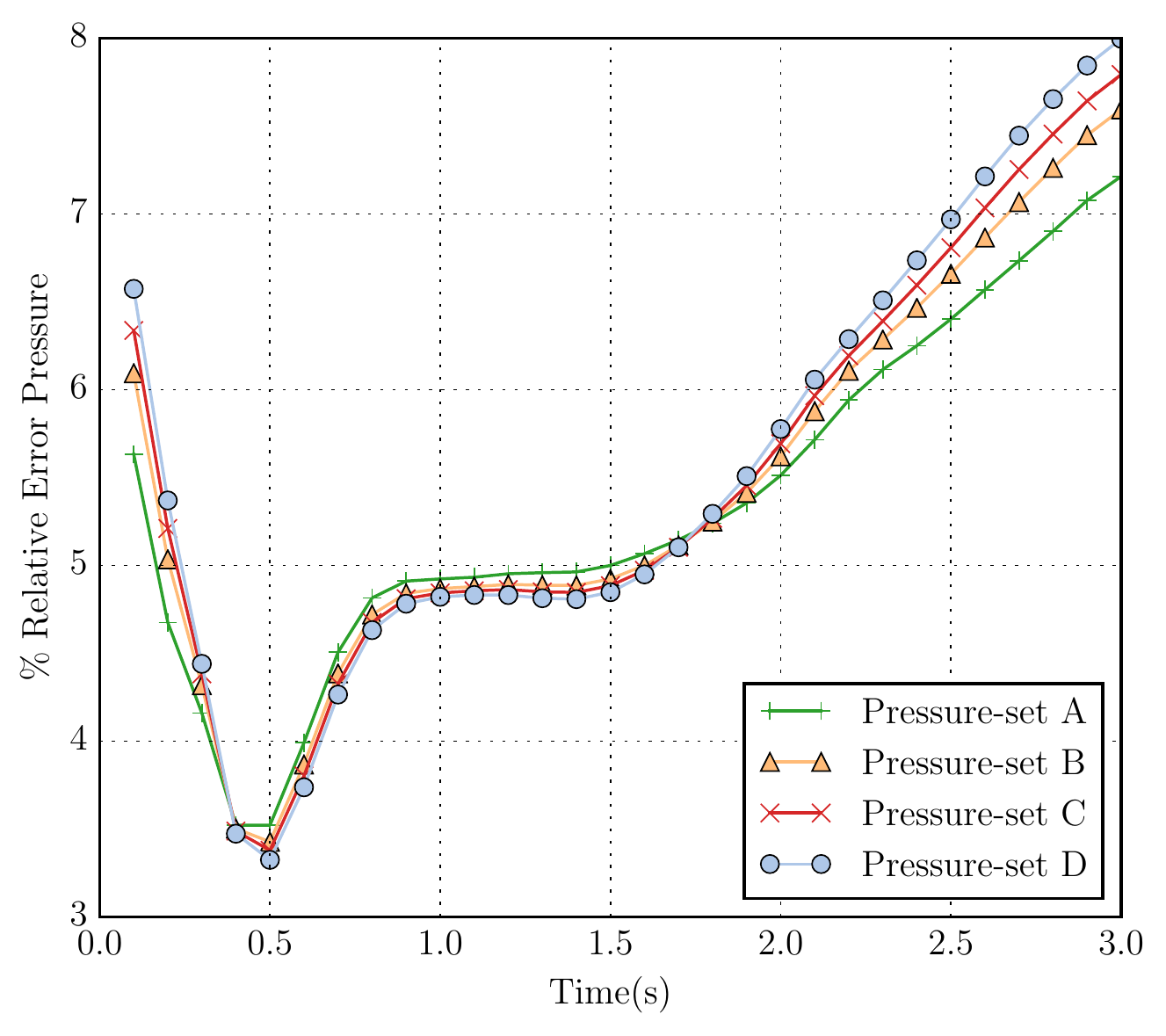}
        \caption{$\% \epsilon_{L^2}(t)$ error of pressure field for the four test sets (\ref{tab:test_points}).}\label{fig:errorl2pres}
  \end{minipage}
\end{figure*}

\begin{figure*}[!tbp]
  \centering 
  \begin{minipage}[b]{0.45\textwidth}
    \includegraphics[width=\textwidth]{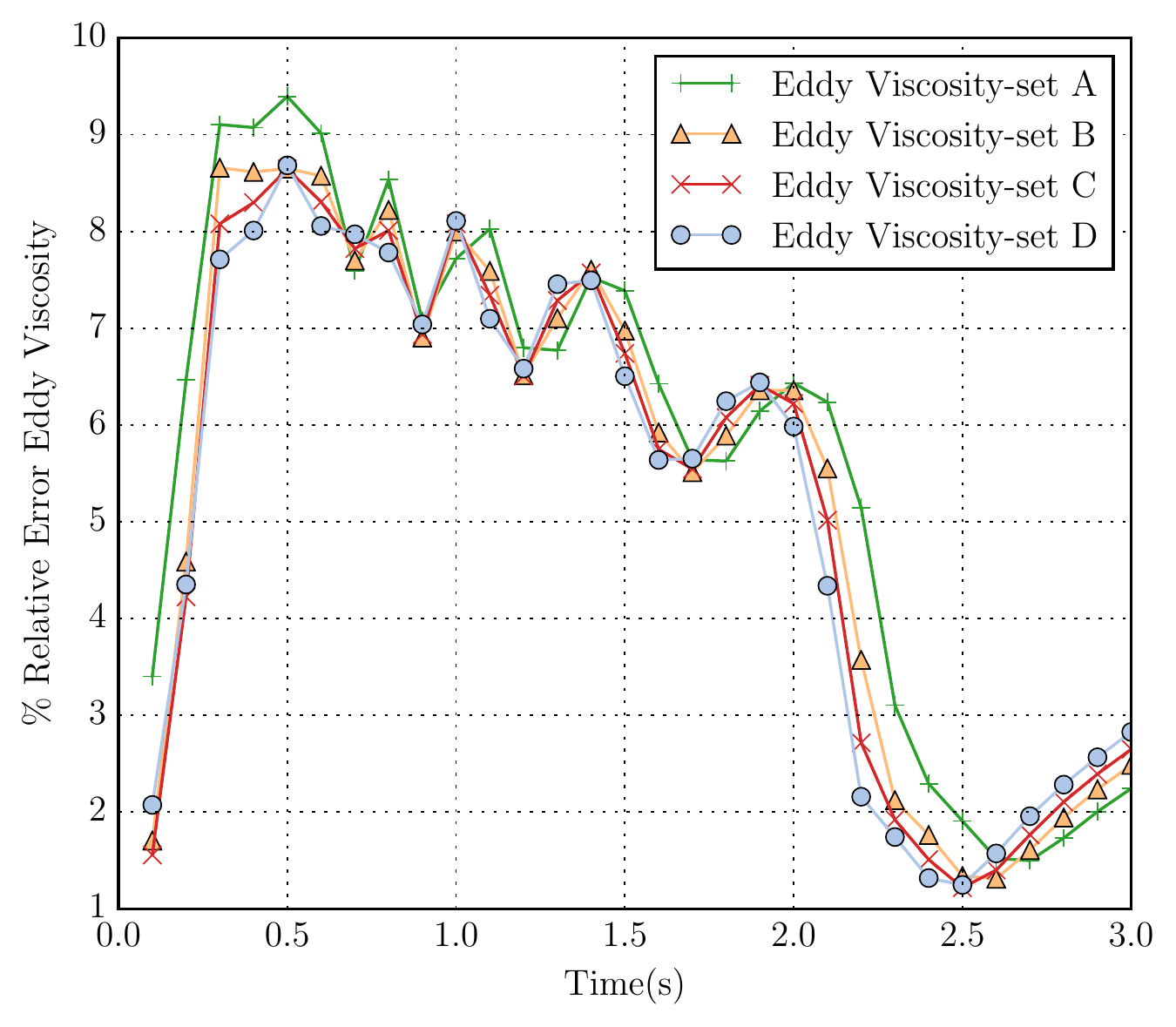}
        \caption{$\% \epsilon_{L^2}(t)$ error of eddy viscosity field for the four test sets (\ref{tab:test_points}).}\label{fig:errorl2eddy}
  \end{minipage}
  \hfill
  \begin{minipage}[b]{0.45\textwidth}  
    \includegraphics[width=\textwidth]{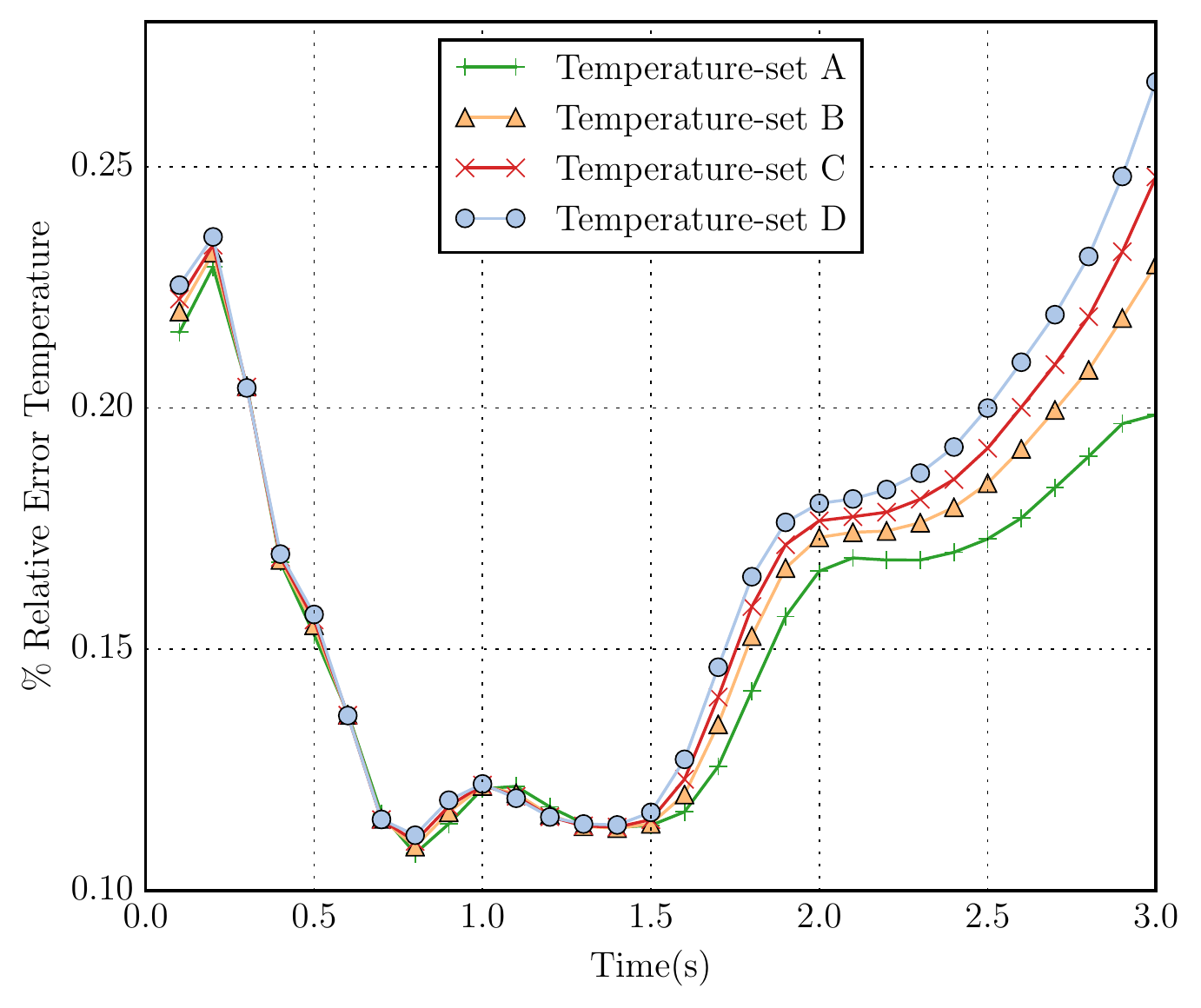}
        \caption{$\% \epsilon_{L^2}(t)$ error of teperature field for the four test sets (\ref{tab:test_points}).}\label{fig:errorl2temp}
  \end{minipage}
\end{figure*}

The ROM is tested on four different sets of velocity inlet values that belong to the training space (table \ref{tab:test_points}). Figures (\ref{fig:errorl2vel}), (\ref{fig:errorl2pres}), (\ref{fig:errorl2eddy}) and (\ref{fig:errorl2temp}) show the relative $\% \epsilon_{L^2}(t)$ error of the velocity, pressure, eddy viscosity and temperature fields for the four test sets. It can be observed that, for the non-linear velocity field, in total, the ROM performs best on set D while for pressure, temperature and eddy viscosity fields, the approximation results in almost similar errors. It can be also observed from figure (\ref{fig:errorl2vel}) that the relative error on the non-linear velocity field grows for sets that lie closer to the lower limit of the training space. As the training space is being enriched with new training points, the ROM appears to perform better, in terms of the velocity field, and therefore set D exhibits the lowest relative error.

\begin{table*}[tbp]  % title of Table
\centering % used for centering table
\tabcolsep=0.09cm
\begin{tabular}{c |c} % centered columns 
\hline\hline %inserts double horizontal lines
$\bm{U_m}$ (m/s) & $\bm{U_b}$ (m/s) \\ [0.5ex] % inserts table
%heading
\hline % inserts single horizontal line
(0.535,0,0) & (0,0,-0.715)  \\ % inserting body of the table
(0.545,0,0) & (0,0,-0.725)  \\
(0.555,0,0) & (0,0,-0.735)  \\
(0.565,0,0) & (0,0,-0.745)  \\
(0.575,0,0) & (0,0,-0.755)  \\
(0.585,0,0) & (0,0,-0.765)  \\
(0.595,0,0) & (0,0,-0.775)  \\
(0.605,0,0) & (0,0,-0.785)  \\
(0.615,0,0) & (0,0,-0.795)  \\
(0.625,0,0) & (0,0,-0.805)  \\% [1ex] adds vertical  
\hline %inserts single line
\end{tabular}
\caption{Sampling points for the parameters.}
\label{tab:sampling_points} % is used to refer this table in the text
\end{table*}

\begin{table*}[tbp]  % title of Table
\centering % used for centering table
\tabcolsep=0.09cm
\begin{tabular}{c |c |c |c |c} % centered columns 
\hline\hline %inserts double horizontal lines
set & $\bm{U_m}$ (m/s) & $\bm{U_b}$ (m/s) & $Re_m$ & $Re_b$ \\ [0.5ex] % inserts table
%heading
\hline % inserts single horizontal line
A & (0.550,0,0) & (0,0,-0.730) & 77000 & 58400   \\ % inserting body of the table
B & (0.570,0,0) & (0,0,-0.750) & 79800 & 60000\\
C & (0.580,0,0) & (0,0,-0.760) & 81200 & 60800\\
D & (0.590,0,0) & (0,0,-0.770) & 82600 & 61600 \\
% [1ex] adds vertical space
\hline %inserts single line
\end{tabular}
\caption{Testing points for the parameters.}
\label{tab:test_points} % is used to refer this table in the text
\end{table*}

\begin{table*}[!tbp]  % title of Table
\centering % used for centering table
\tabcolsep=0.09cm
\begin{tabular}{c | c | c} % centered columns 
\hline\hline %inserts double horizontal lines
& FOM & ROM \\ [0.5ex] % inserts table
%heading
\hline % inserts single horizontal line
$\partial/\partial t$ & Backward & Backward  \\ % inserting body of the table
$\nabla\cdot$ & Upwind, Central & Upwind, Central  \\
$\nabla^2$ & Central & Central \\
\hline %inserts single line
\end{tabular}
\caption{Numerical Schemes for FOM and ROM.}
\label{tab:schemes} % is used to refer this table in the text
\end{table*}
\raggedbottom

\setlength\intextsep{0pt}  
\setlength\intextsep{0pt}

% \setlength\intextsep{0pt} 
% \begin{figure*}[!tbp] 
%   \centering 
%   \begin{minipage}[b]{0.40\textwidth}
%     \includegraphics[width=\textwidth]{figures/modes.pdf}
%     \caption{Number of POD modes for velocity, temperature, pressure, supremizer and nut fields needed to reach approximately $ 99.9\%$ of energy as a function of mesh size, using URANS.}\label{fig:ransmodes}
%   \end{minipage}
%   \hfill
%   \begin{minipage}[b]{0.40\textwidth}
%     \includegraphics[width=\textwidth]{figures/modesles.pdf}
%     \caption{Number of POD modes for velocity, temperature, pressure, supremizer and nut fields needed to reach approximately $99.9\%$ of energy as a function of mesh size, using LES.}\label{fig:lesmodes}
%   \end{minipage} 
% \end{figure*}
% \setlength\intextsep{0pt} 

\subsection{Nested and Standard POD methods}\label{sec:num_exp}
 
In this subsection, a comparison between nested POD and standard POD methods for the test set D (table \ref{tab:test_points}) is presented. The boundary conditions for the FOM are shown in table (\ref{fig:boundary_con}). In the standard POD, the POD is applied directly on the global snapshots matrices as explained in \ref{sec:training}.

\begin{figure*}[!tbp]
  \centering 
  \begin{minipage}[b]{0.40\textwidth}
    \includegraphics[width=\textwidth]{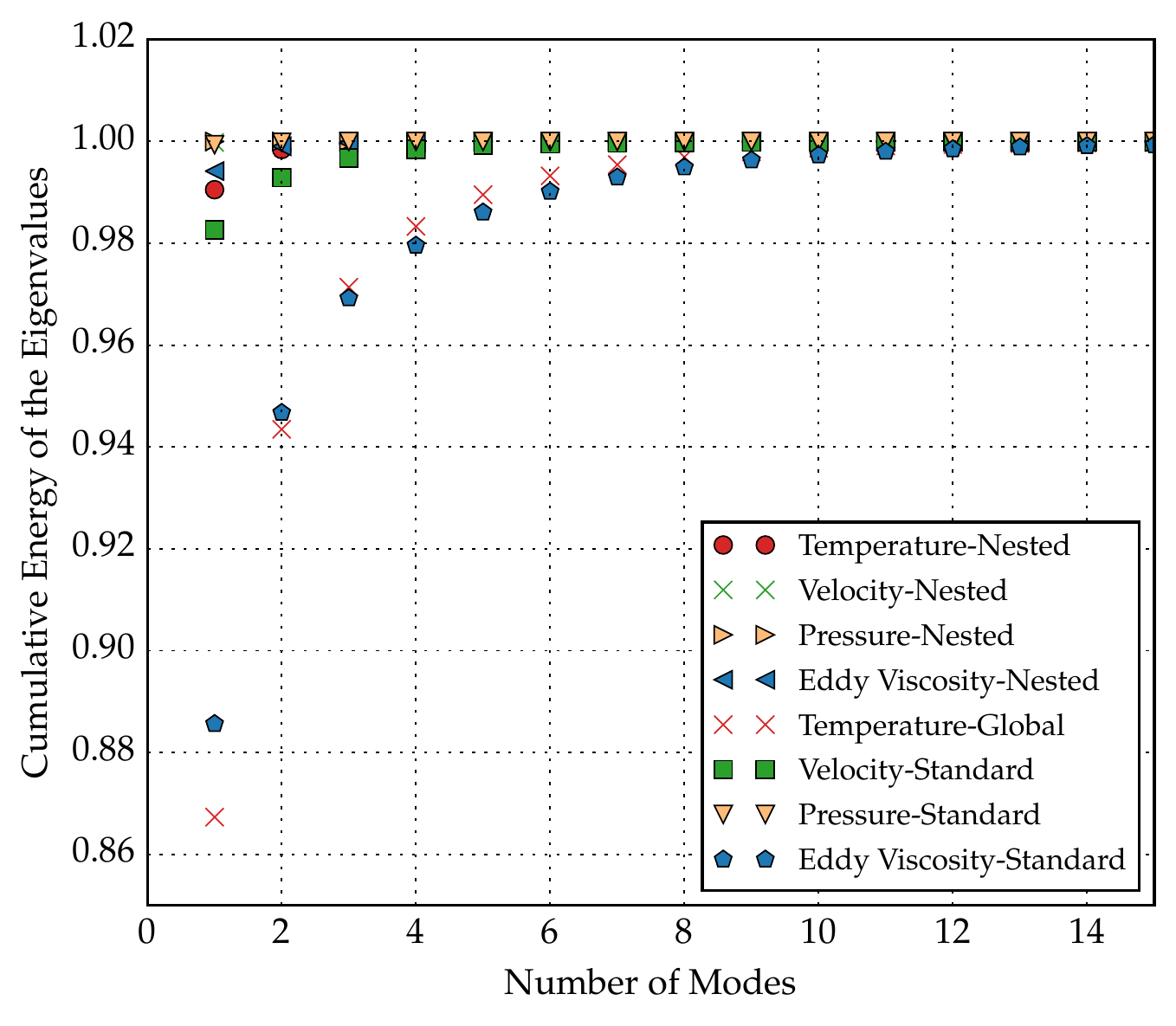}
    \caption{Cumulative energy of the eigenvalues for temperature, velocity, pressure and eddy viscosity fields for nested and standard POD methods, respectively.}\label{fig:eigenvalues_case2}\label{fig:energyeigen}
  \end{minipage}
  \hfill
  \begin{minipage}[b]{0.40\textwidth}  
    \includegraphics[width=\textwidth]{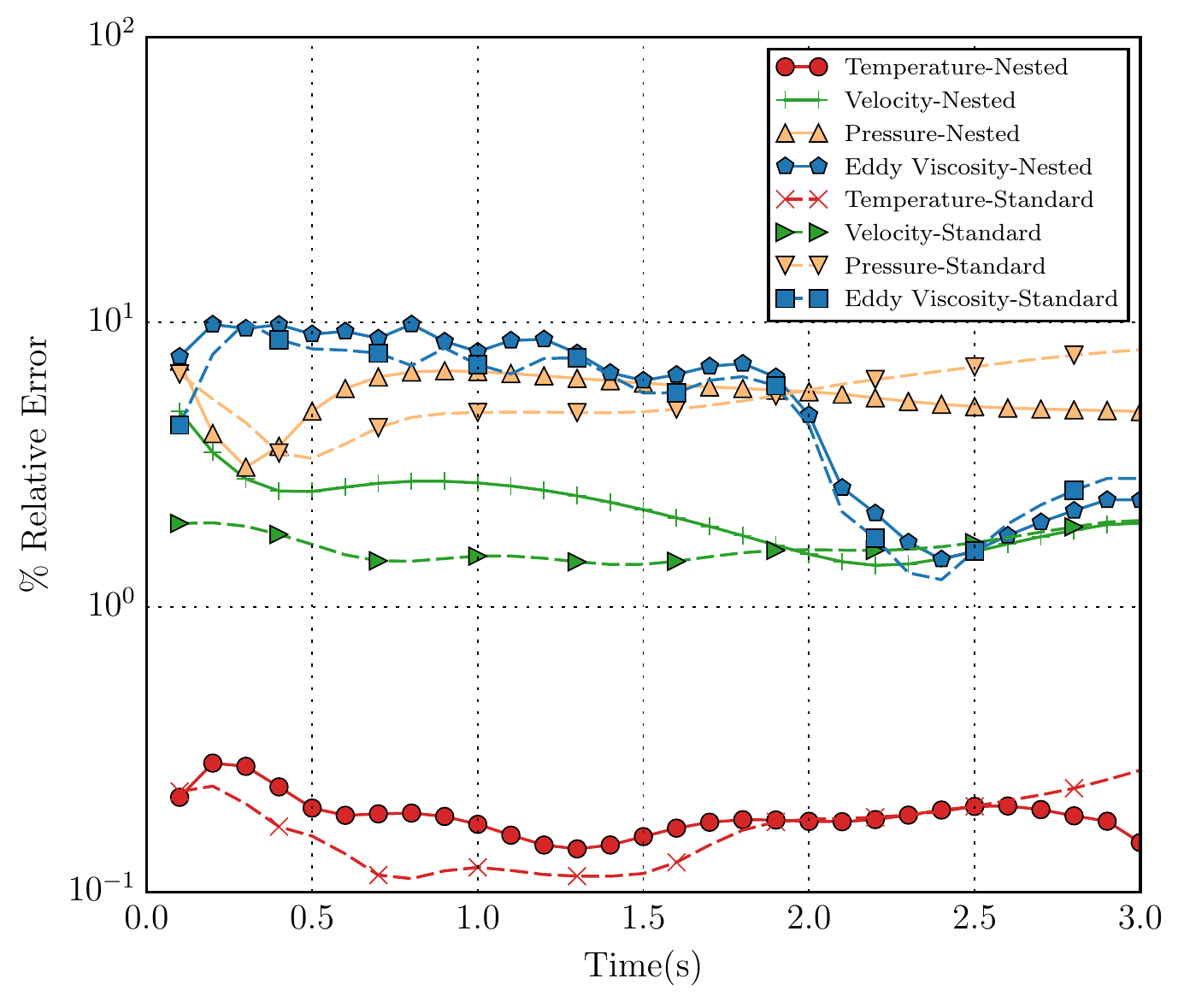}
    \caption{ $\% \epsilon_{L^2}(t)$ error for temperature, velocity, pressure and eddy viscosity fields for the nested and standard POD methods, respectively.}\label{fig:errorl2}
  \end{minipage}
\end{figure*}

In the nested POD method, $10$ local snapshot matrices are constructed (\ref{nestedmatrix}), one for each sampling pair, and the POD is applied on each of them individually. The resulting bases functions are then truncated using the energy criterion (\ref{eq:energycr}) in order to retain approximately $99.9\%$ of the total energy and the remained bases are weighted by their eigenvalues. According to this criterion, out of the $30$ bases functions, only $10$ are retained and weighted for each sampling pair, giving a reduced dimension for the final global matrix of $N_{u}^{rn}=10\cdot10=100$, according to (\ref{smallglobal}). Finally, the POD is applied on the weighted global matrix and the same procedure as in the standard POD is followed. The final decay of the eigenvalues can be seen in figure (\ref{fig:energyeigen}). For the purposes of the comparison, the same number of modes for both the standard POD and the nested POD have been chosen. This makes $N_u^r=6$, $N_p^r=N_{sup}^{r}=10$, $N_{\theta}^r=11$ and $N_{\nu_t}^r=10$. The numerical effort of the nested POD can be calculated as $\mathscr{O}([N_{t}^3\cdot N_p +[N_{u}^{\mbox{ nested }}]^3) =\mathscr{O}(30^3\cdot 10 + [10\cdot10]^3)\approx 2\cdot 10^6$, while for the standard POD $\mathscr{O}([N_{t}^3\cdot N_p +[N_{u}^{\mbox{ nested }}]^3) =\mathscr{O}([30\cdot10]^3)\approx 2\cdot 10^7$. Therefore, the nested POD is one order of magnitude faster. For cases with a higher number of snapshots and sampling points, the smaller numerical effort of the nested POD is more apparent, see for example the results in \cite{doi:10.1002/gamm.201610011}.   

In the online phase, the ROM is evaluated in a new set of values (test set D) $U_m=0.59$ m/s and $U_b=0.77$ m/s. The ROM is tested for $3$s with timestep $dt = 0.0025$s. The total execution time of the FOM is $13782.1$s on a single processor while of the ROM $11.02$s. This sets the ROM approximately 1259 times faster than the FOM (table \ref{tab:cpu}). 

In figure (\ref{fig:errorl2}), the relative $\%\epsilon_{L^2}$ error between the FOM and the ROM constructed using standard POD and nested POD methods is shown. This is calculated as:

\begin{equation}
\epsilon_{L^2}(t) = \frac{||X_{FOM}(t)-X_{ROM}(t)||_{L^2(\Omega)}}{||X_{FOM}(t)||_{L^2(\Omega)}}\% .
\end{equation} 

According to this figure the relative error for the standard POD method is slightly lower during the first 2 seconds of the simulation followed by similar or slightly better performance of the nested POD for the rest of the simulation. Tables (\ref{tab:errorstatnested}) and (\ref{tab:errorstat}) summarise a few statistics, showing minimum, maximum and average error for nested and standard POD methods, where, the average error appears to be very close for both cases. The pressure and eddy viscosity fields appear to have the largest relative error according to figure (\ref{fig:errorl2}). For pressure, this can be attributed to the lack of an explicit equation, as in the ROM, the pressure bases functions are projected onto the continuity equation using the supremizer method \cite{Rozza2007,Ballarin2014}. An exploitation of a Poisson equation for pressure would probably improve the error \cite{noack2005,Caiazzo2014598,Stabile2017}. As for the eddy viscosity field, the use of non-intrusive, pure data-driven method  could cause the slighlty higher relative error, compared to the velocity and temperature fields, where the exact equations are projected onto the reduced bases. For the temperature field, the energy equation is still weakly coupled with the momentum equation, therefore the error appears smaller for temperature field. 

\begin{figure*}[!tbp]
\begin{minipage}{1\textwidth}
\centering 
\includegraphics[width=0.30\textwidth]{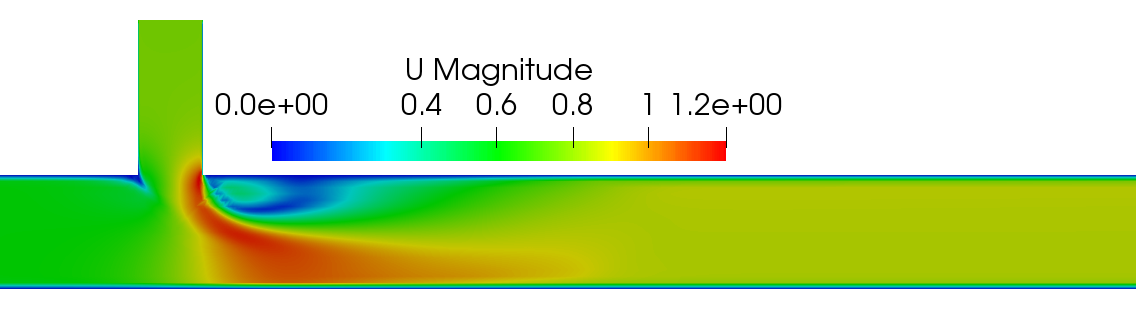}
\includegraphics[width=0.30\textwidth]{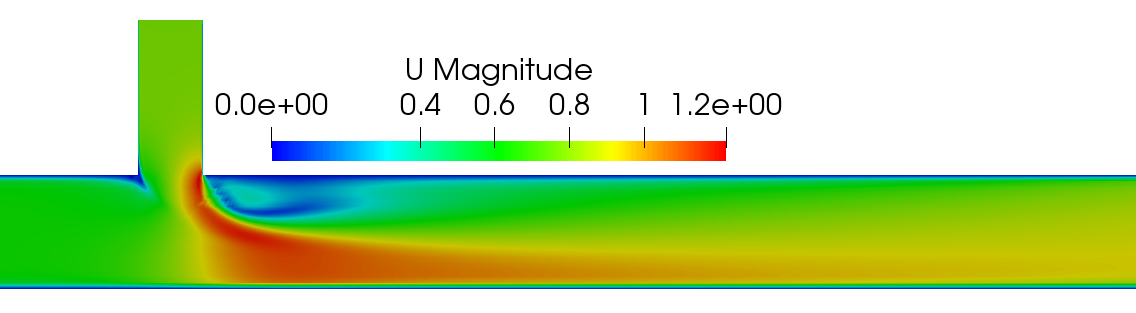}
\includegraphics[width=0.30\textwidth]{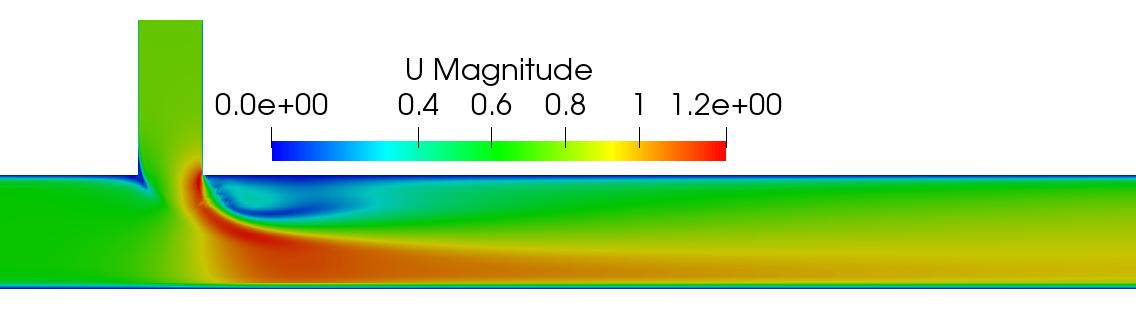}
\includegraphics[width=0.30\textwidth]{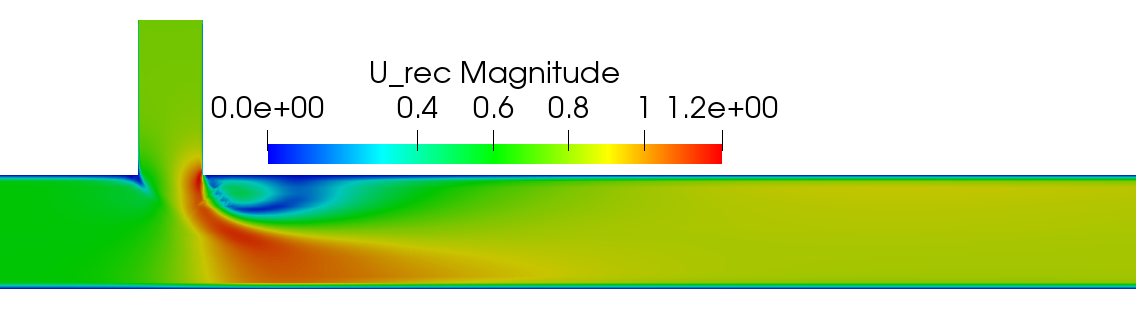}
\includegraphics[width=0.30\textwidth]{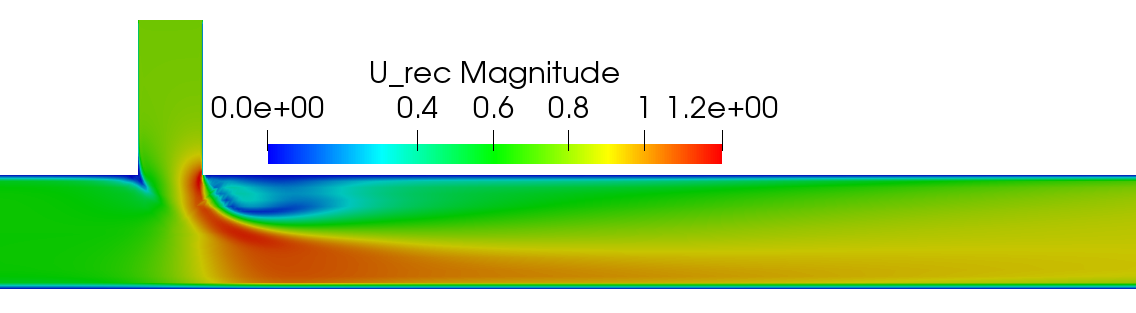}
\includegraphics[width=0.30\textwidth]{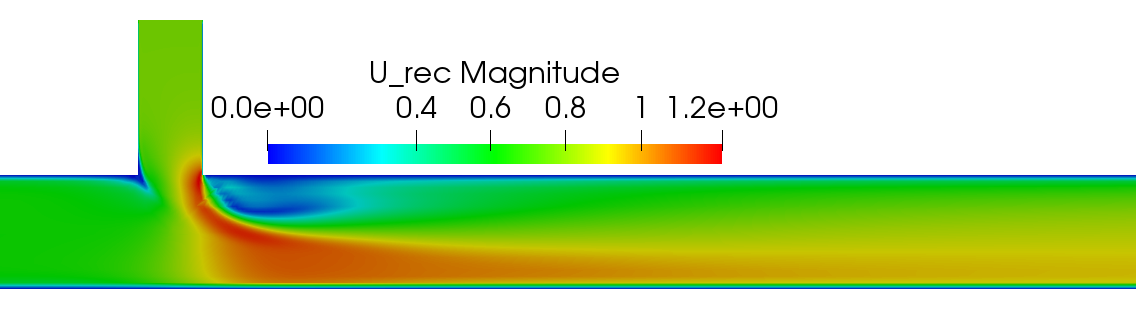}
\includegraphics[width=0.30\textwidth]{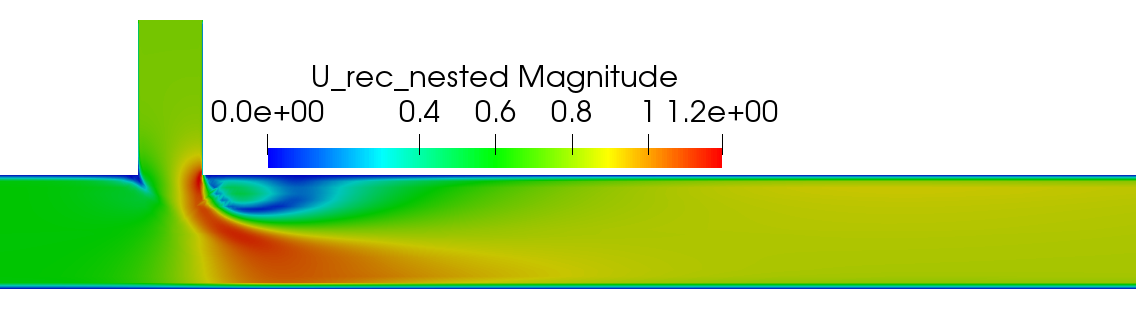}
\includegraphics[width=0.30\textwidth]{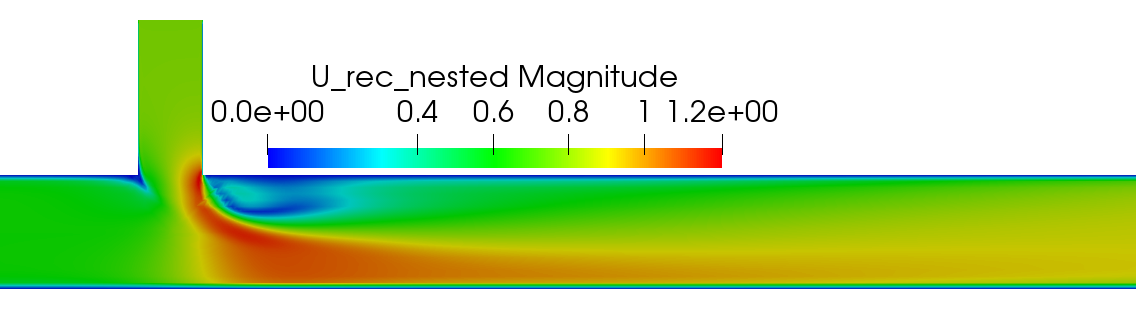}
\includegraphics[width=0.30\textwidth]{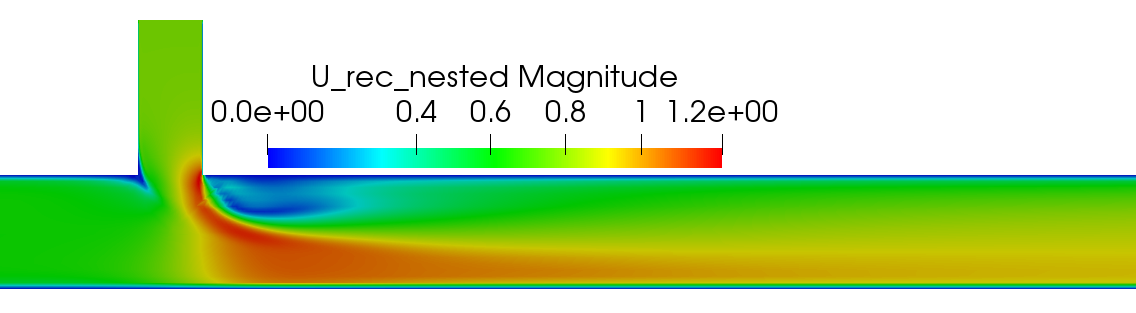}
\
\end{minipage} 
\caption{Comparison of the velocity field of the full order (first row) and reduced order model using standard POD (second row) and nested POD method (third row). The fields are depicted for different time instances equal to $t=0.5 \si{s}, 1.5\si{s}$ and $3 \si{s}$ and increasing from left to right.}\label{fig:comparison_case2_vel}
\end{figure*}

\begin{figure*}[!tbp] 
\begin{minipage}{1\textwidth}
\centering 
\includegraphics[width=0.30\textwidth]{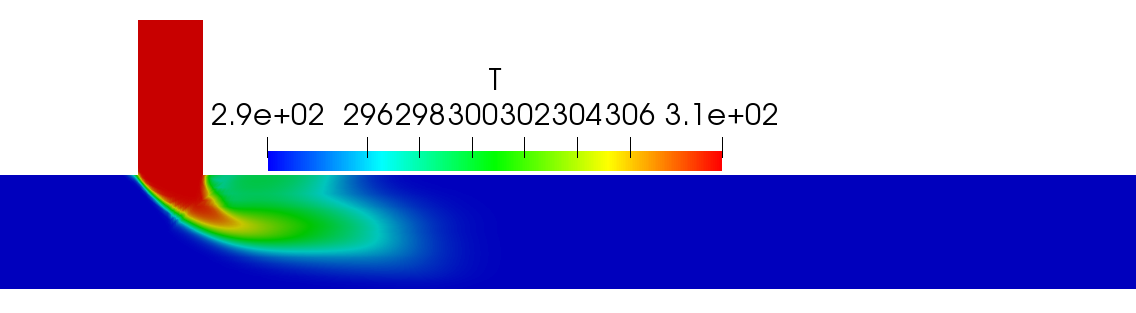}
\includegraphics[width=0.30\textwidth]{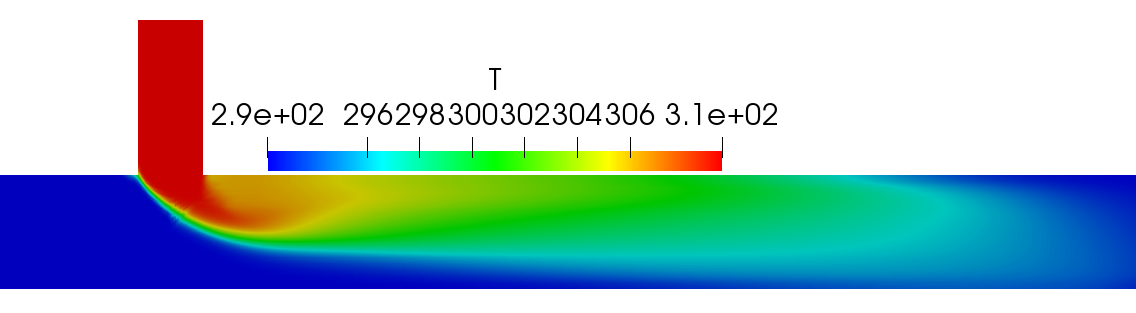}  
\includegraphics[width=0.30\textwidth]{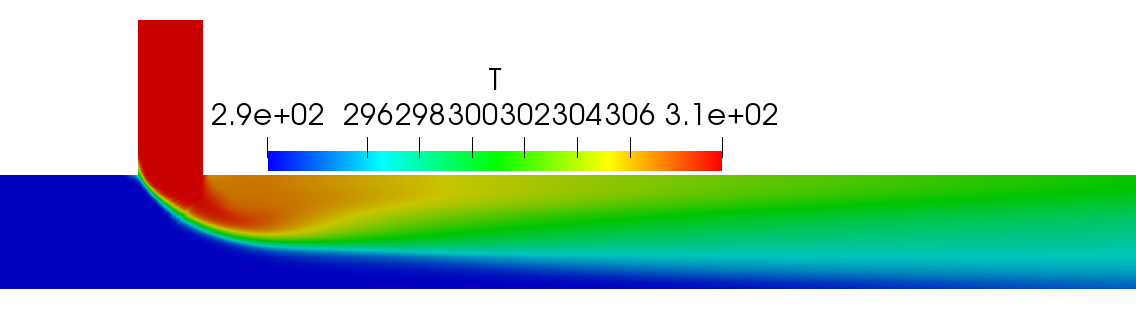} 
\includegraphics[width=0.30\textwidth]{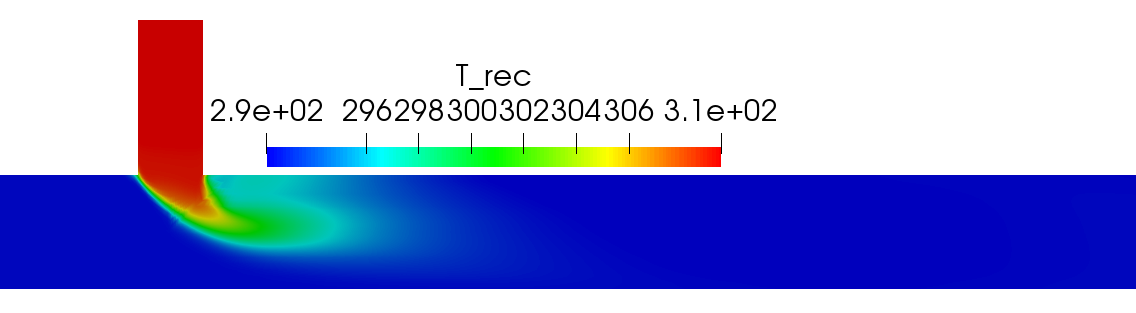}
\includegraphics[width=0.30\textwidth]{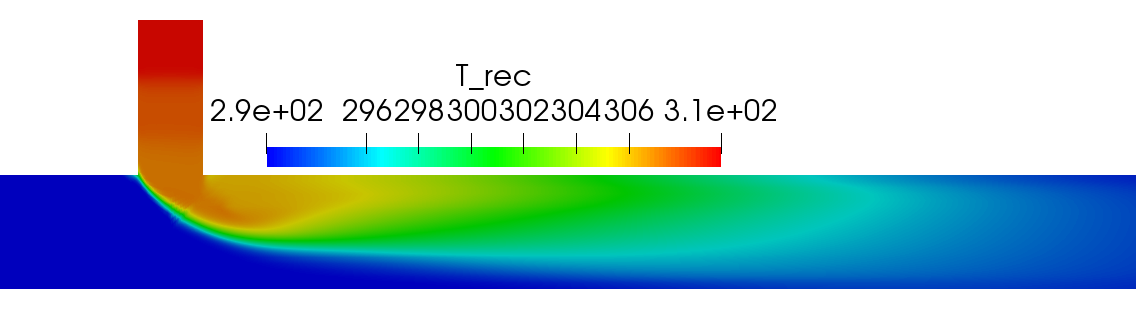}
\includegraphics[width=0.30\textwidth]{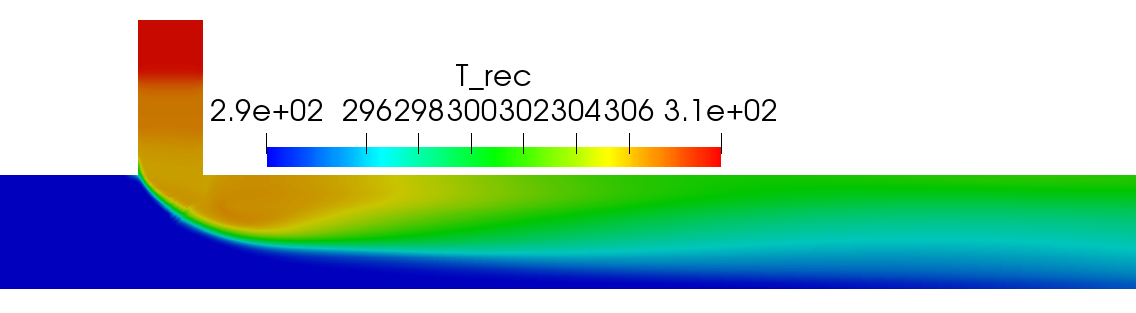}
\includegraphics[width=0.30\textwidth]{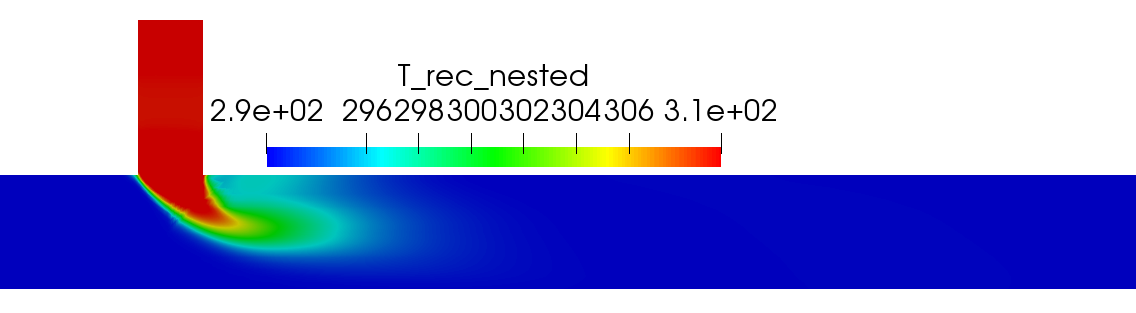}
\includegraphics[width=0.30\textwidth]{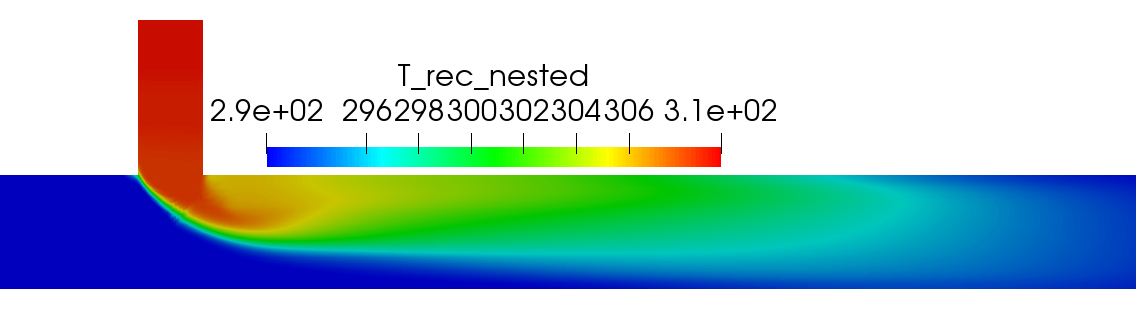}
\includegraphics[width=0.30\textwidth]{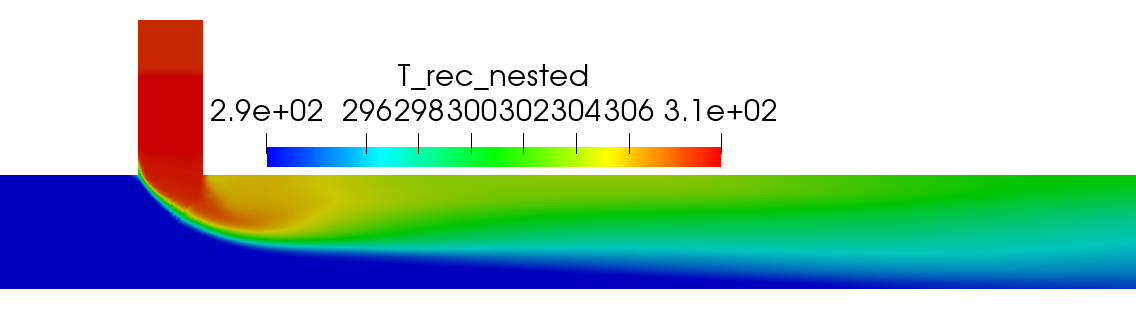}
\
\end{minipage} 
\caption{Comparison of the temperature field of the full order (first row) and reduced order model using standard POD (second row) and nested POD method (third row). The fields are depicted for different time instances equal to $t=0.5 \si{s}, 1.5\si{s}$ and $3 \si{s}$ and increasing from left to right.}\label{fig:comparison_case2_temp}
\end{figure*}

\begin{figure*}[!tbp]
\begin{minipage}{1\textwidth}
\centering 
\includegraphics[width=0.30\textwidth]{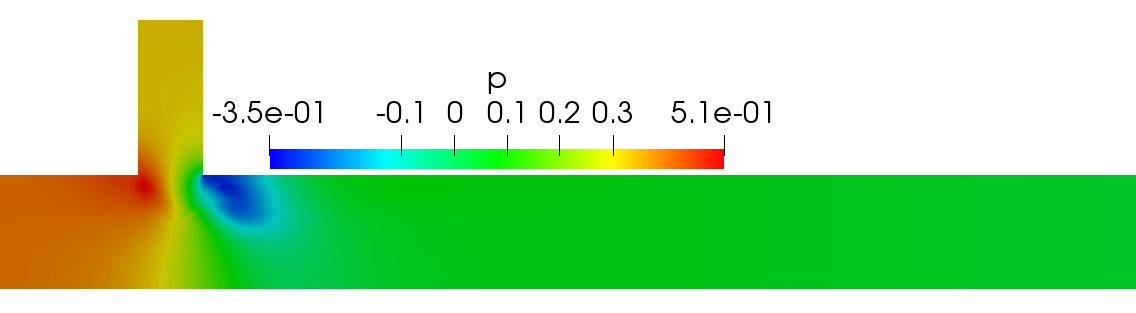}
\includegraphics[width=0.30\textwidth]{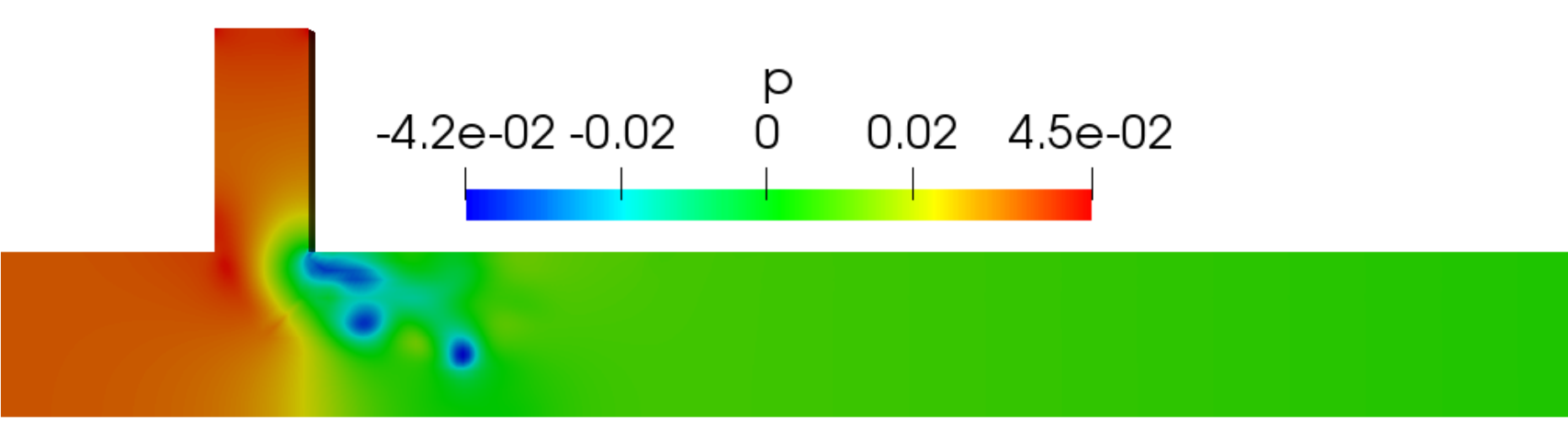}
\includegraphics[width=0.30\textwidth]{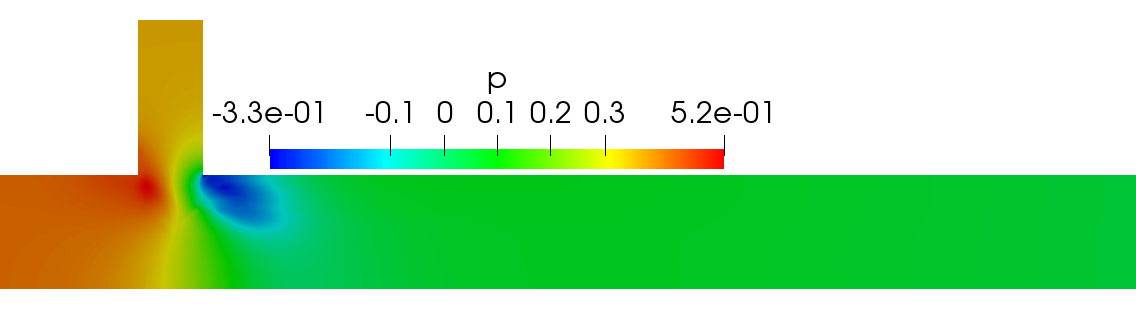}
\includegraphics[width=0.30\textwidth]{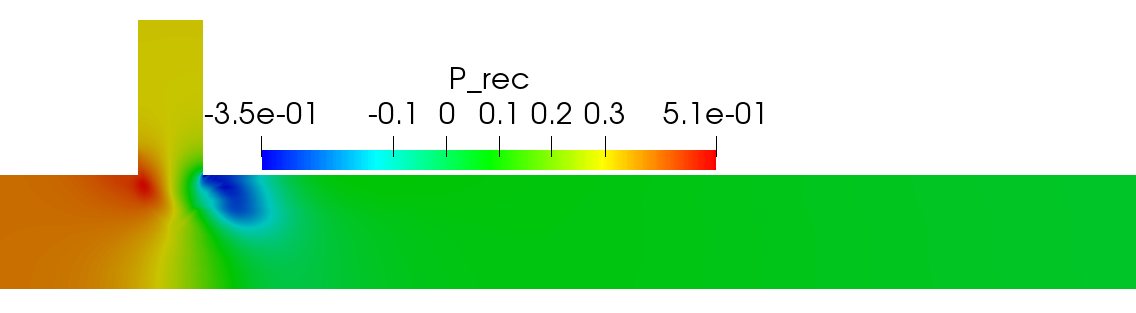}
\includegraphics[width=0.30\textwidth]{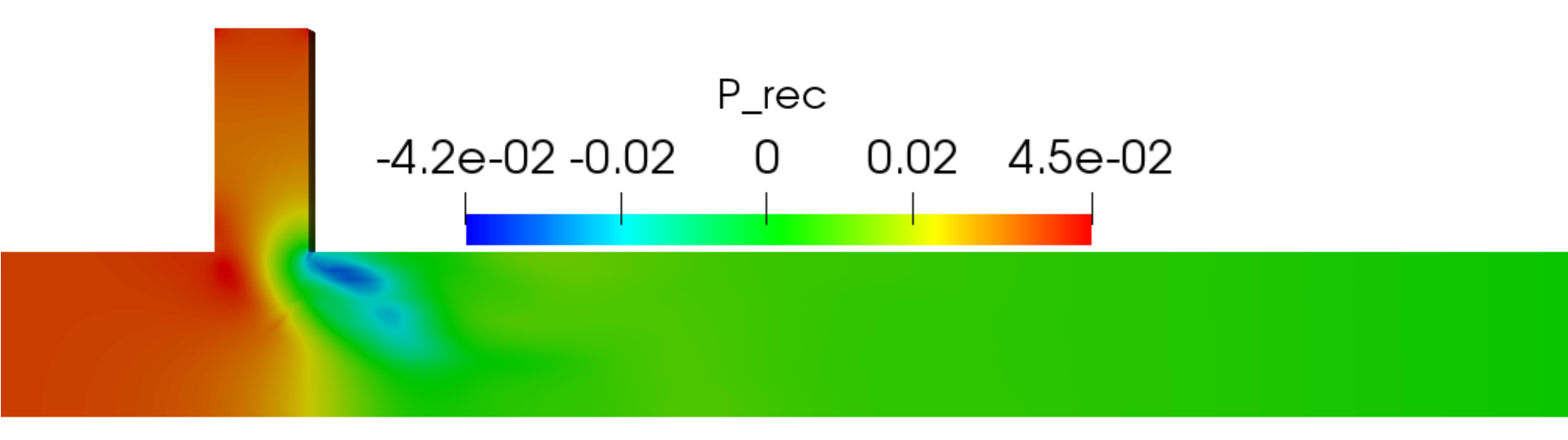}
\includegraphics[width=0.30\textwidth]{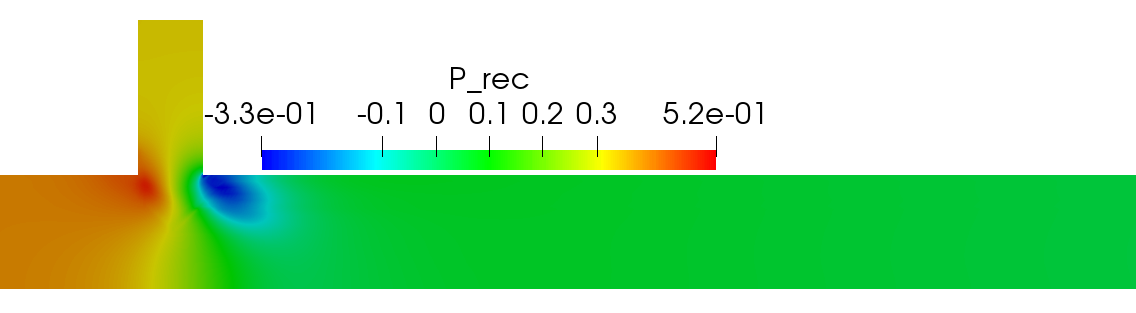}
\includegraphics[width=0.30\textwidth]{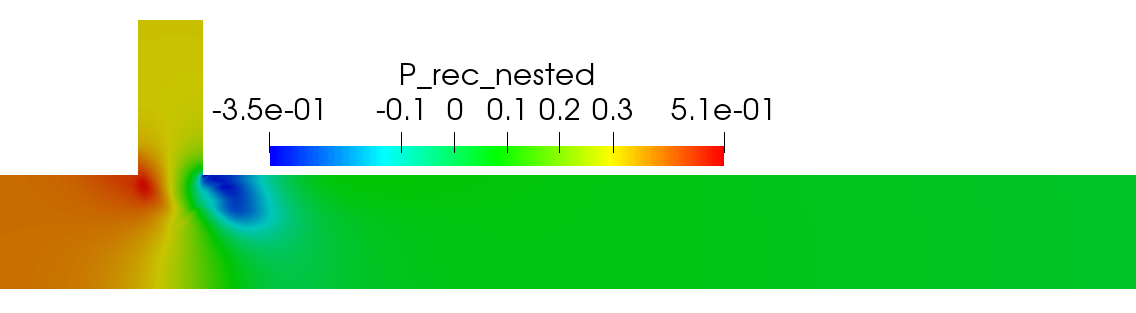}
\includegraphics[width=0.30\textwidth]{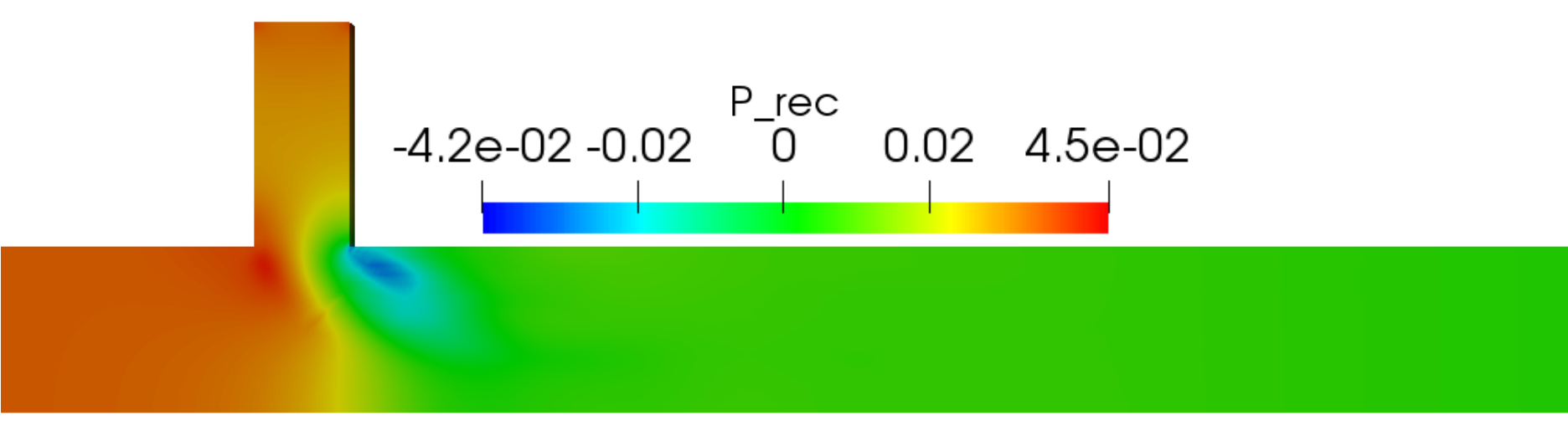}
\includegraphics[width=0.30\textwidth]{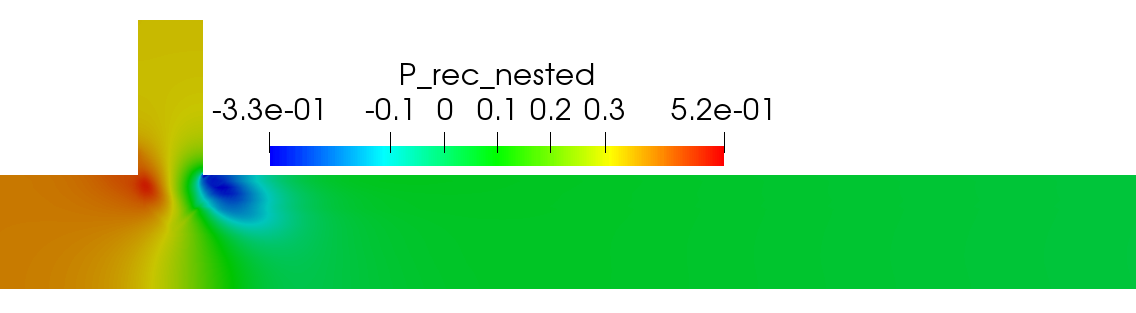}
\
\end{minipage} 
\caption{Comparison of the pressure field of the full order (first row) and reduced order model using standard POD (second row) and nested POD method (third row). The fields are depicted for different time instances equal to $t=0.5 \si{s}, 1.5\si{s}$ and $3 \si{s}$ and increasing from left to right.}\label{fig:comparison_case2_press}
\end{figure*}

\begin{figure*}[!tbp]
\begin{minipage}{1\textwidth}
\centering 
\includegraphics[width=0.30\textwidth]{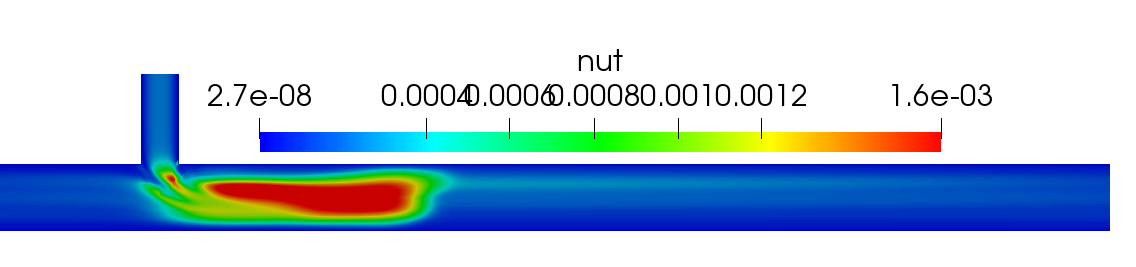}
\includegraphics[width=0.30\textwidth]{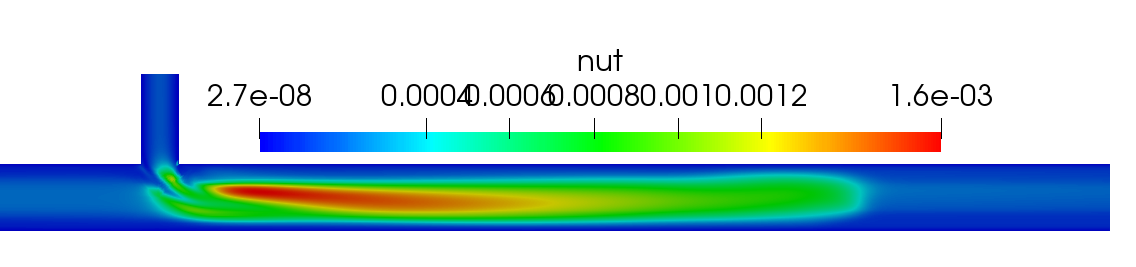}
\includegraphics[width=0.30\textwidth]{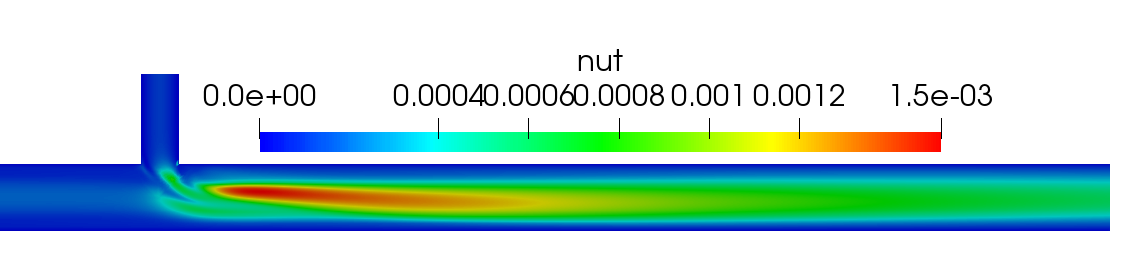}
\includegraphics[width=0.30\textwidth]{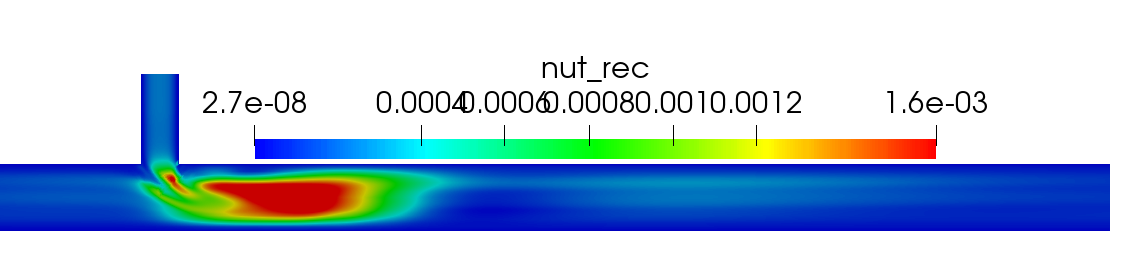}
\includegraphics[width=0.30\textwidth]{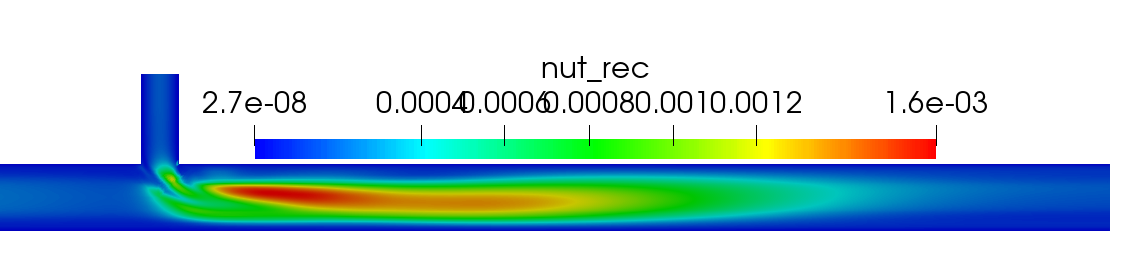}
\includegraphics[width=0.30\textwidth]{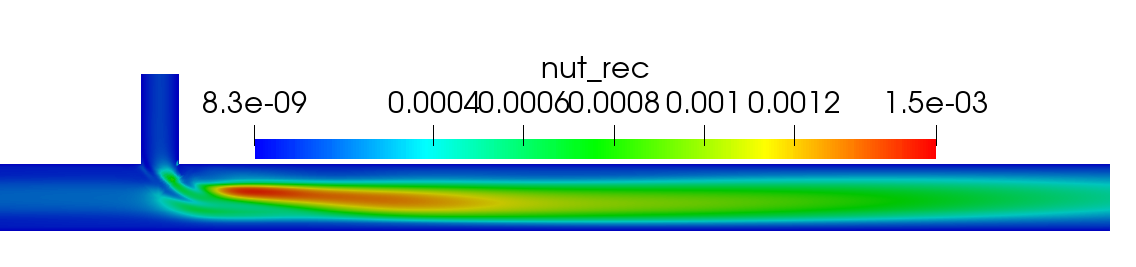}
\includegraphics[width=0.30\textwidth]{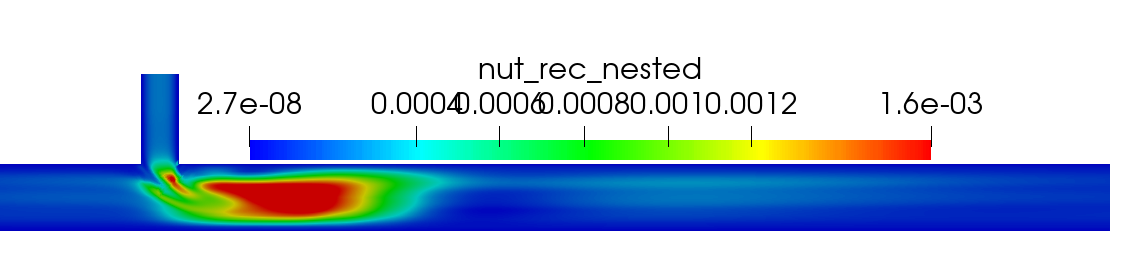}
\includegraphics[width=0.30\textwidth]{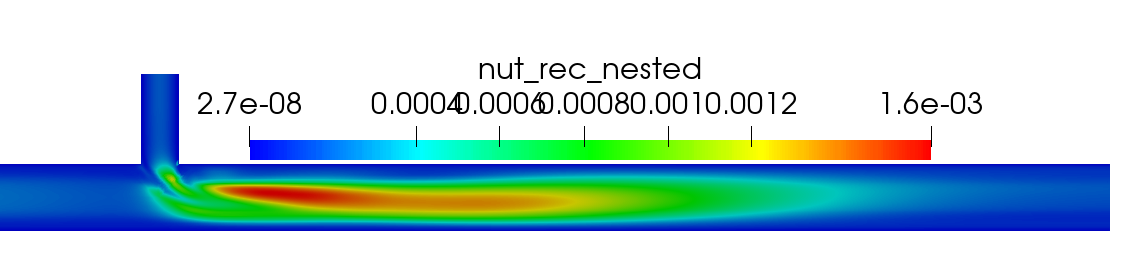}
\includegraphics[width=0.30\textwidth]{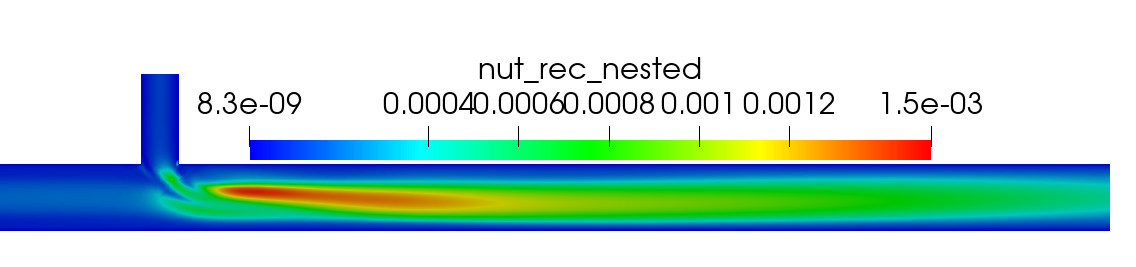}
\
\end{minipage} 
\caption{ Comparison of the eddy viscosity field of the full order (first row), reduced order model with standard POD (second row) and reduced order model with nested POD (third row). The fields are depicted for different time instances equal to $t=0.5 \si{s}, 1.5\si{s}$ and $3 \si{s}$ and increasing from left to right.}\label{fig:comparison_case2_nut}
\end{figure*}

\begin{figure*}[!tbp]
\begin{minipage}{1\textwidth}    
\centering 
\includegraphics[width=0.30\textwidth]{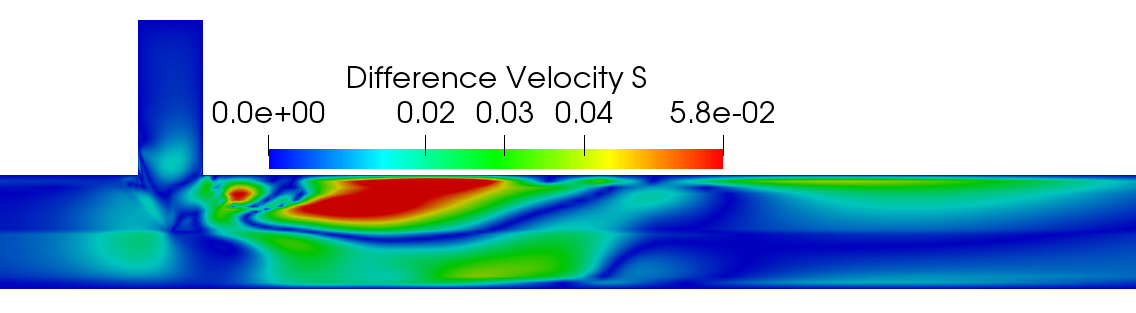}
\includegraphics[width=0.30\textwidth]{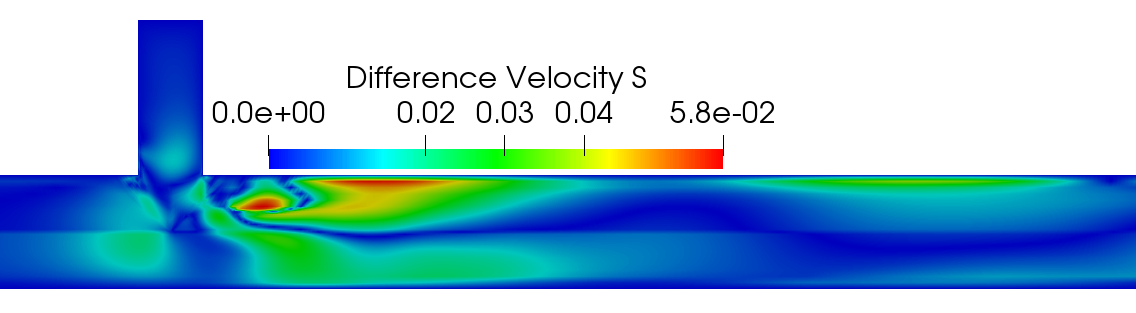}
\includegraphics[width=0.30\textwidth]{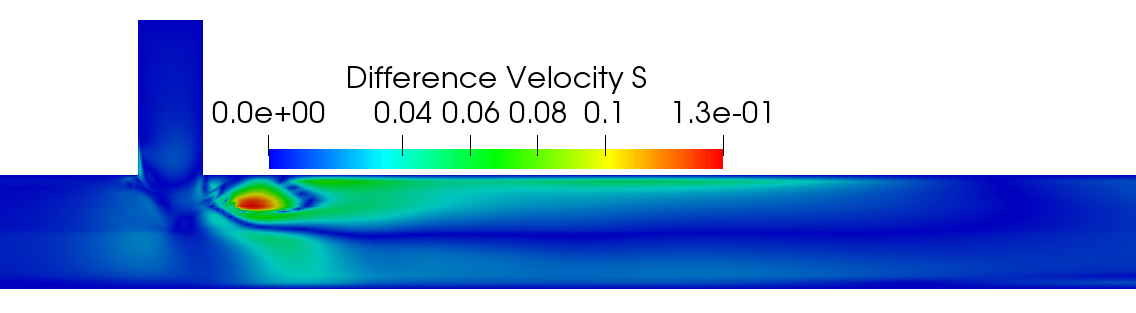}
\includegraphics[width=0.30\textwidth]{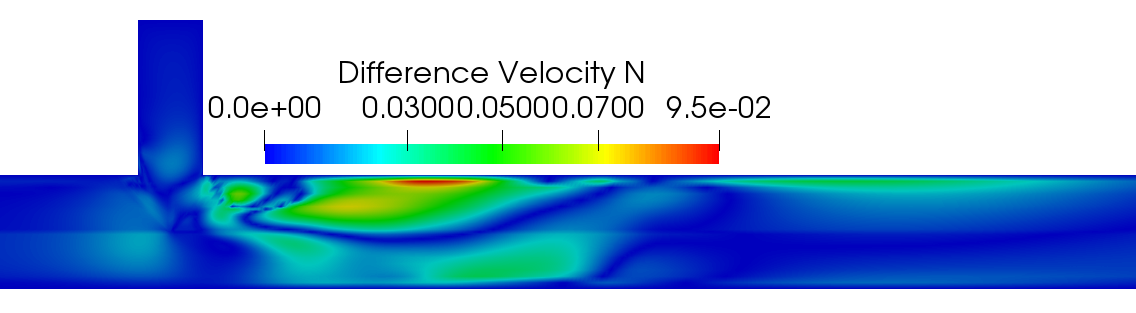}
\includegraphics[width=0.30\textwidth]{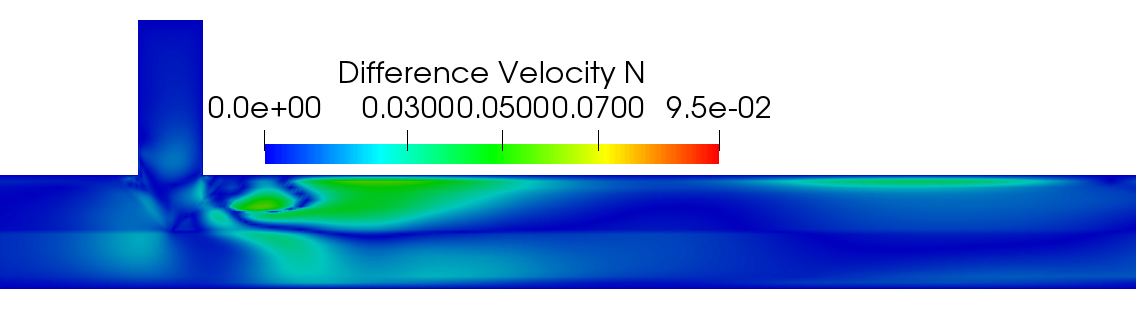}
\includegraphics[width=0.30\textwidth]{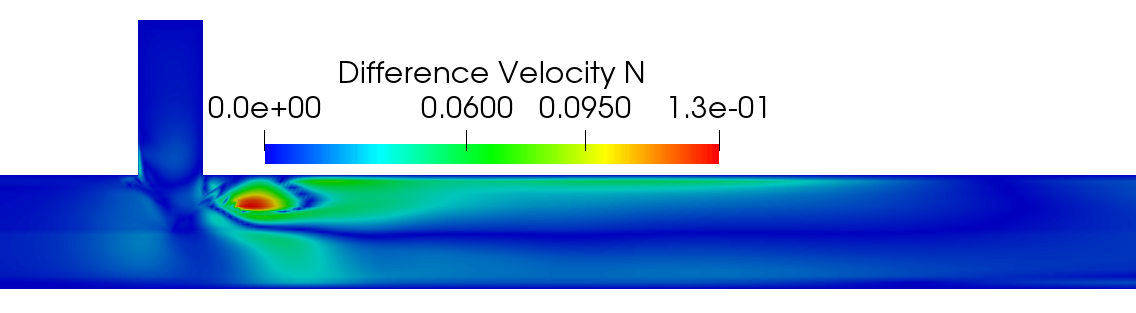}
\includegraphics[width=0.30\textwidth]{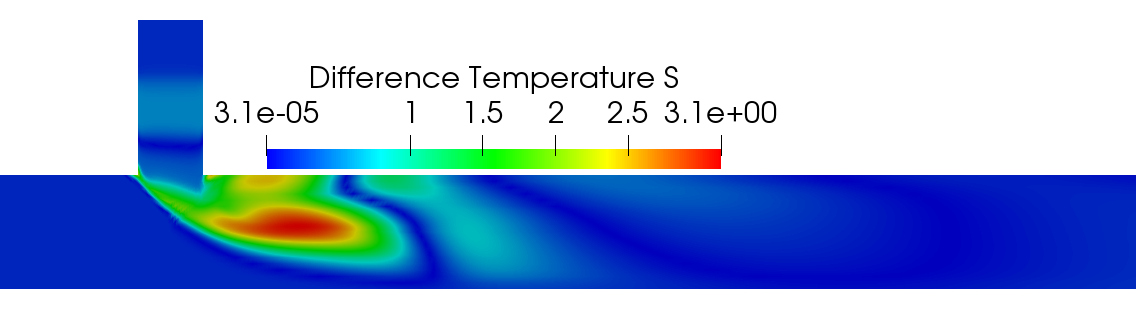}
\includegraphics[width=0.30\textwidth]{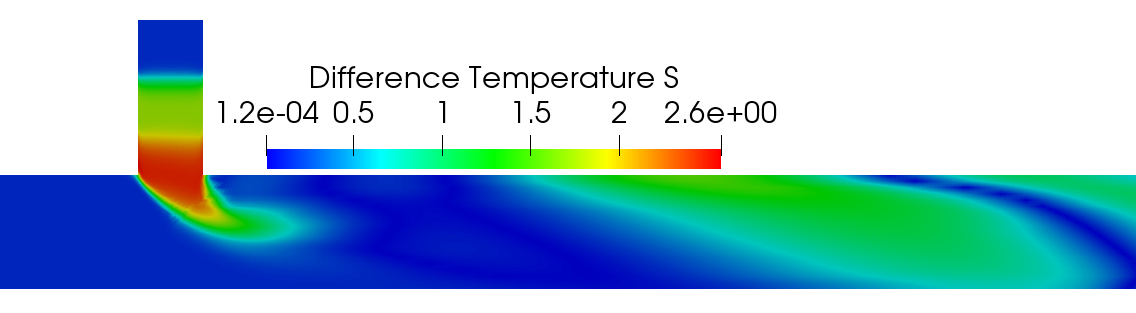}
\includegraphics[width=0.30\textwidth]{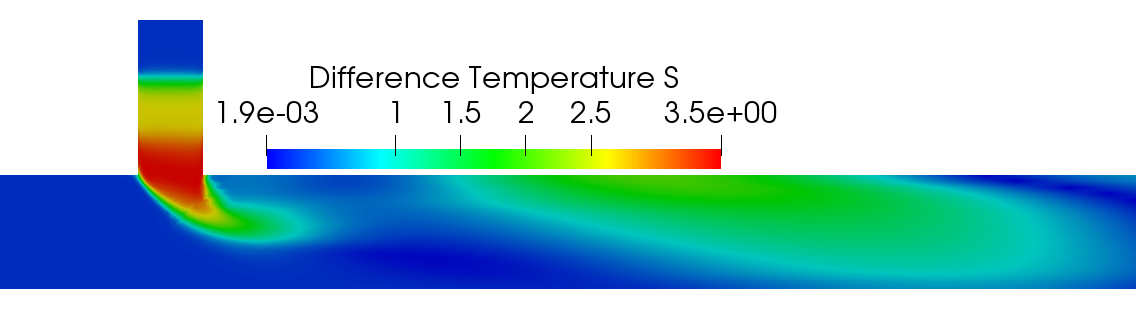}
\includegraphics[width=0.30\textwidth]{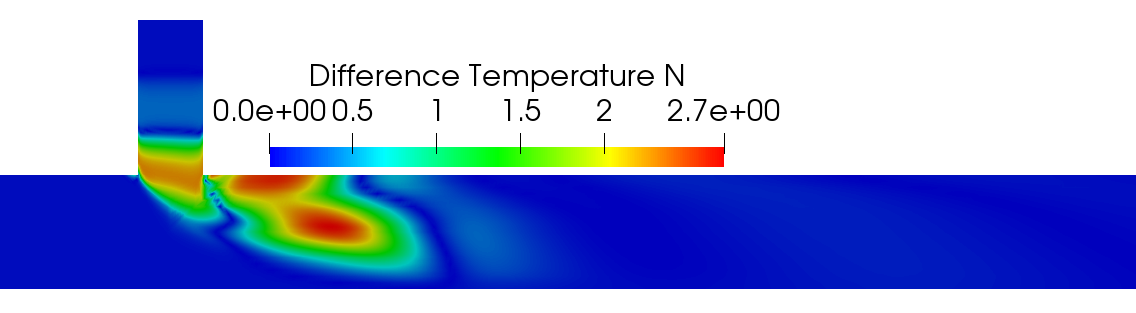}
\includegraphics[width=0.30\textwidth]{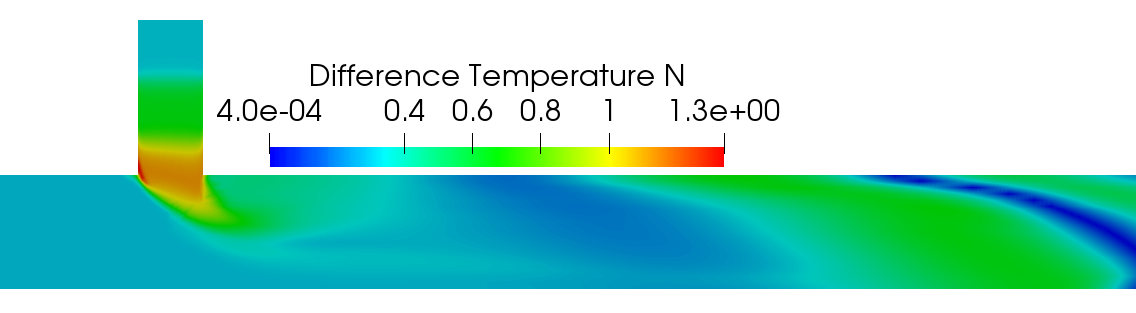}
\includegraphics[width=0.30\textwidth]{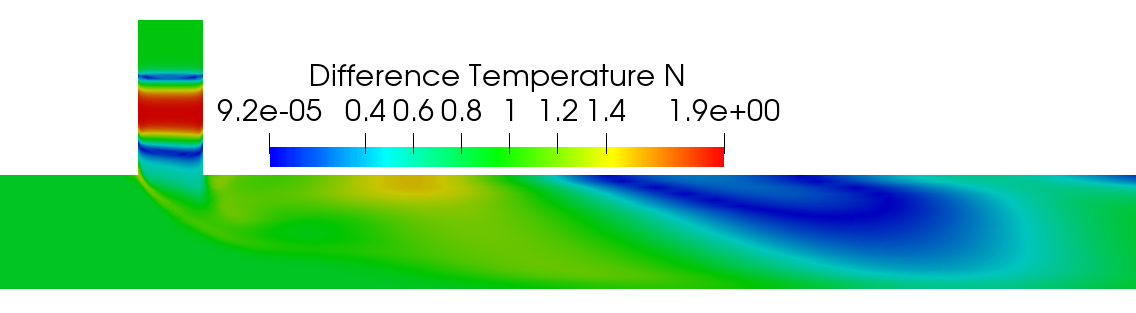}
\includegraphics[width=0.30\textwidth]{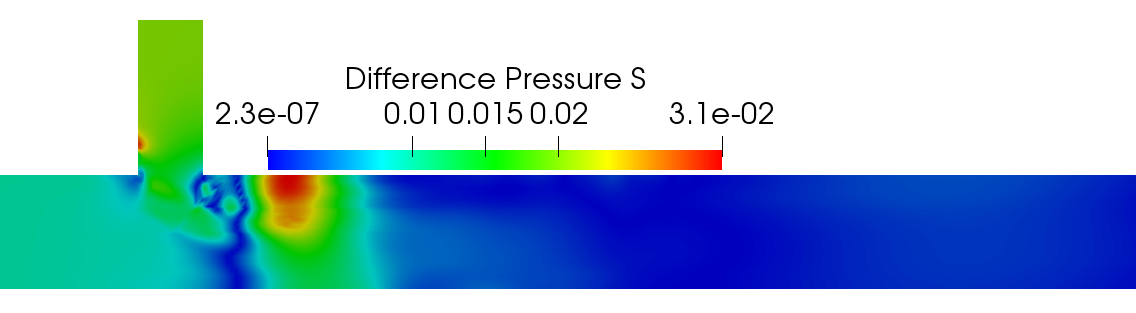}
\includegraphics[width=0.30\textwidth]{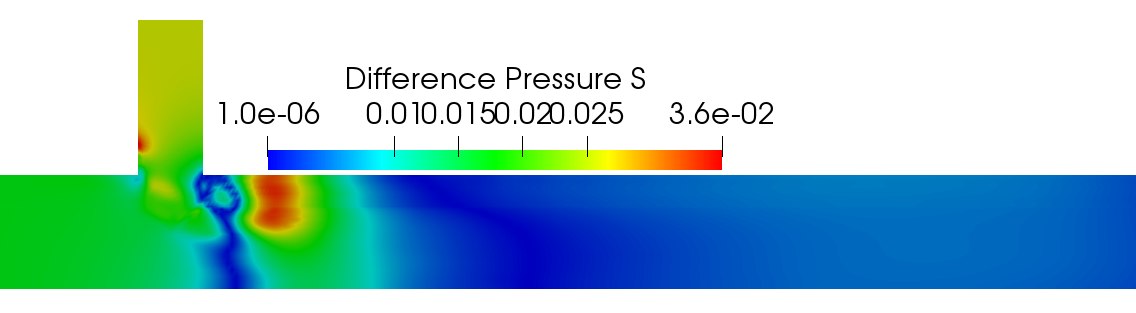}
\includegraphics[width=0.30\textwidth]{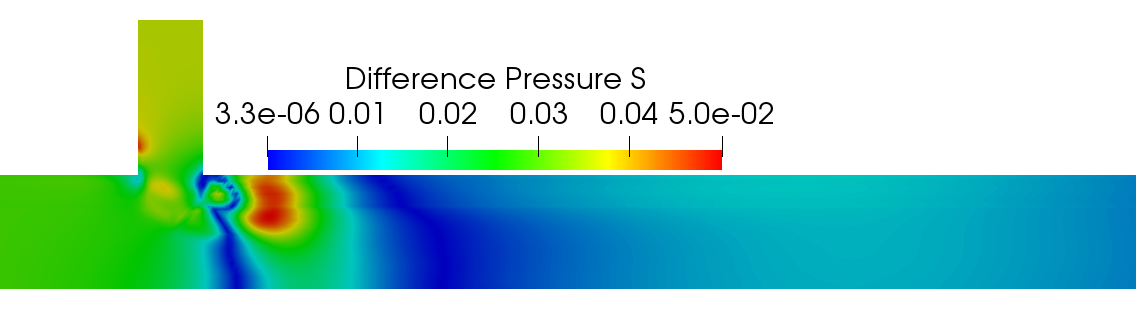}
\includegraphics[width=0.30\textwidth]{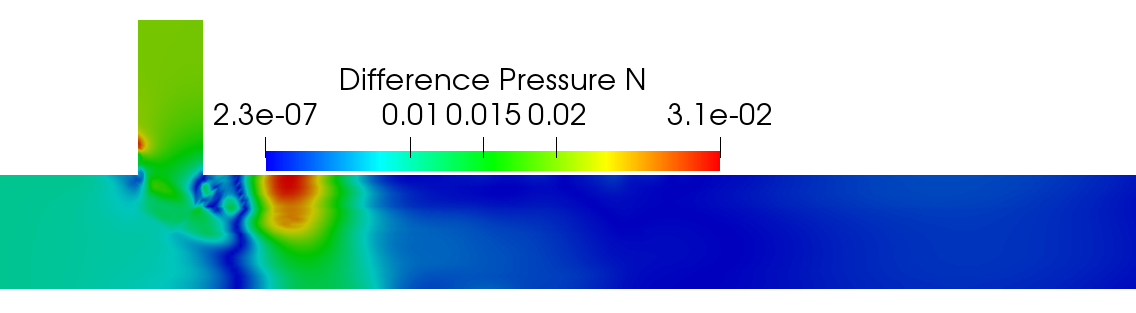}
\includegraphics[width=0.30\textwidth]{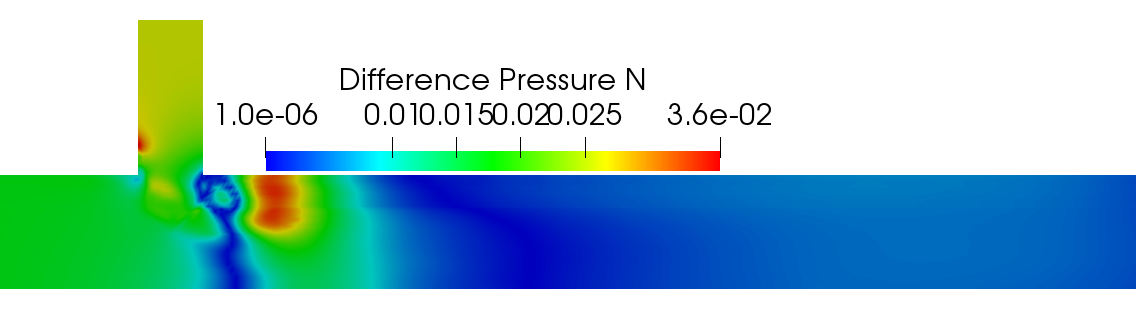}
\includegraphics[width=0.30\textwidth]{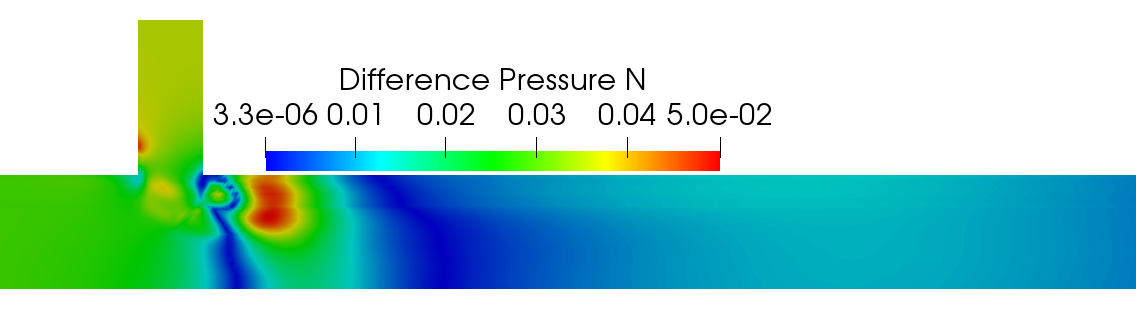}
\
\end{minipage} 
\caption{Difference between the full order and reduced order standard POD (odd rows) and nested POD (even rows) for velocity, temperature and pressure fields, respectively. The fields are depicted for different time instances equal to $t=0.5 \si{s}, 1.5\si{s}$ and $3 \si{s}$ and increasing from left to right.}\label{fig:comparison_case2_diff}
\end{figure*}

\begin{figure*}[!tbp]
\begin{minipage}{1\textwidth}
\centering 
\includegraphics[width=0.30\textwidth]{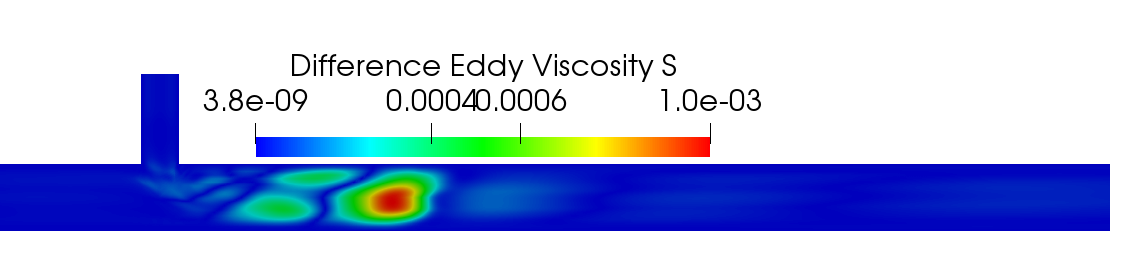}
\includegraphics[width=0.30\textwidth]{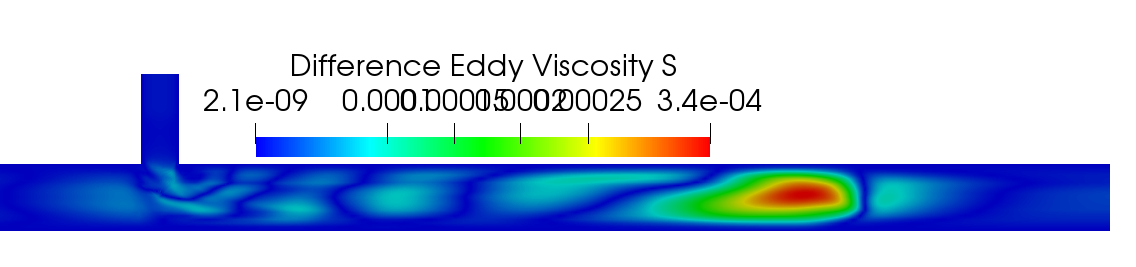}
\includegraphics[width=0.30\textwidth]{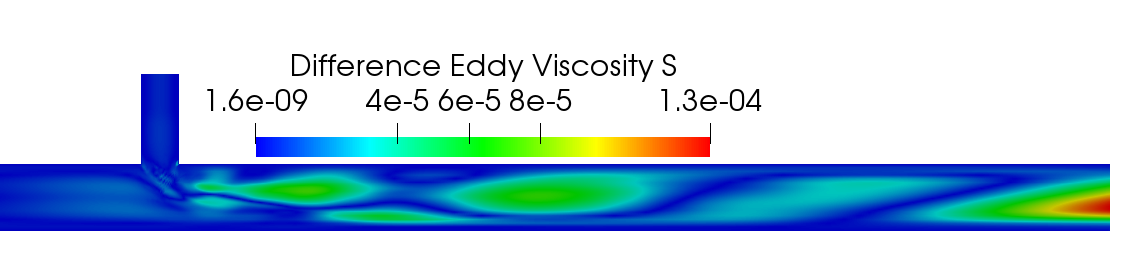}
\includegraphics[width=0.30\textwidth]{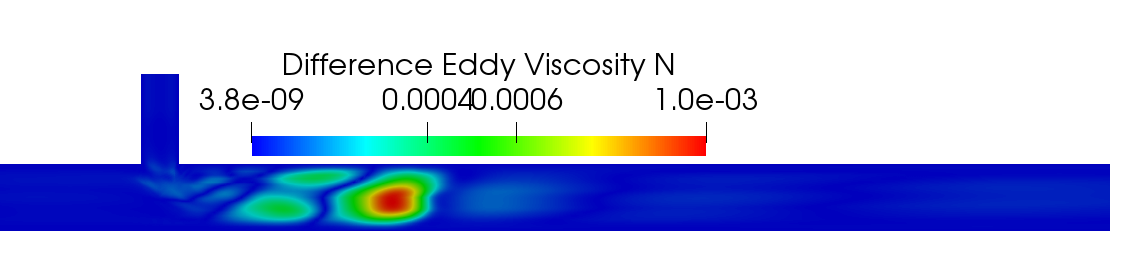}
\includegraphics[width=0.30\textwidth]{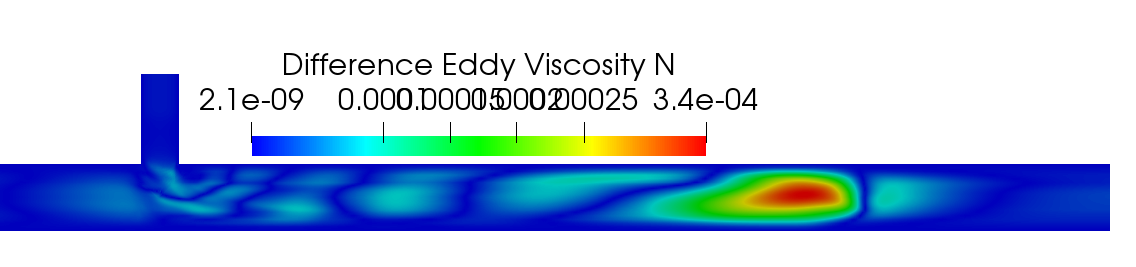}
\includegraphics[width=0.30\textwidth]{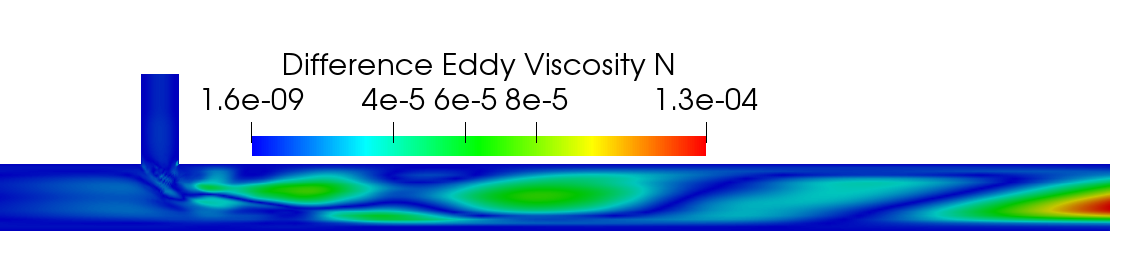}
\
\end{minipage} 
\caption{Difference between the full order and reduced order standard POD (first row) and nested POD (second row) for the eddy viscosity field. The fields are depicted for different time instances equal to $t=0.5 \si{s}, 1.5\si{s}$ and $3 \si{s}$ and increasing from left to right.}\label{fig:comparison_case2_diff_nut}
\end{figure*}

\begin{figure*}[!tbp]
\begin{minipage}{1\textwidth}     
\centering 
\includegraphics[width=0.30\textwidth]{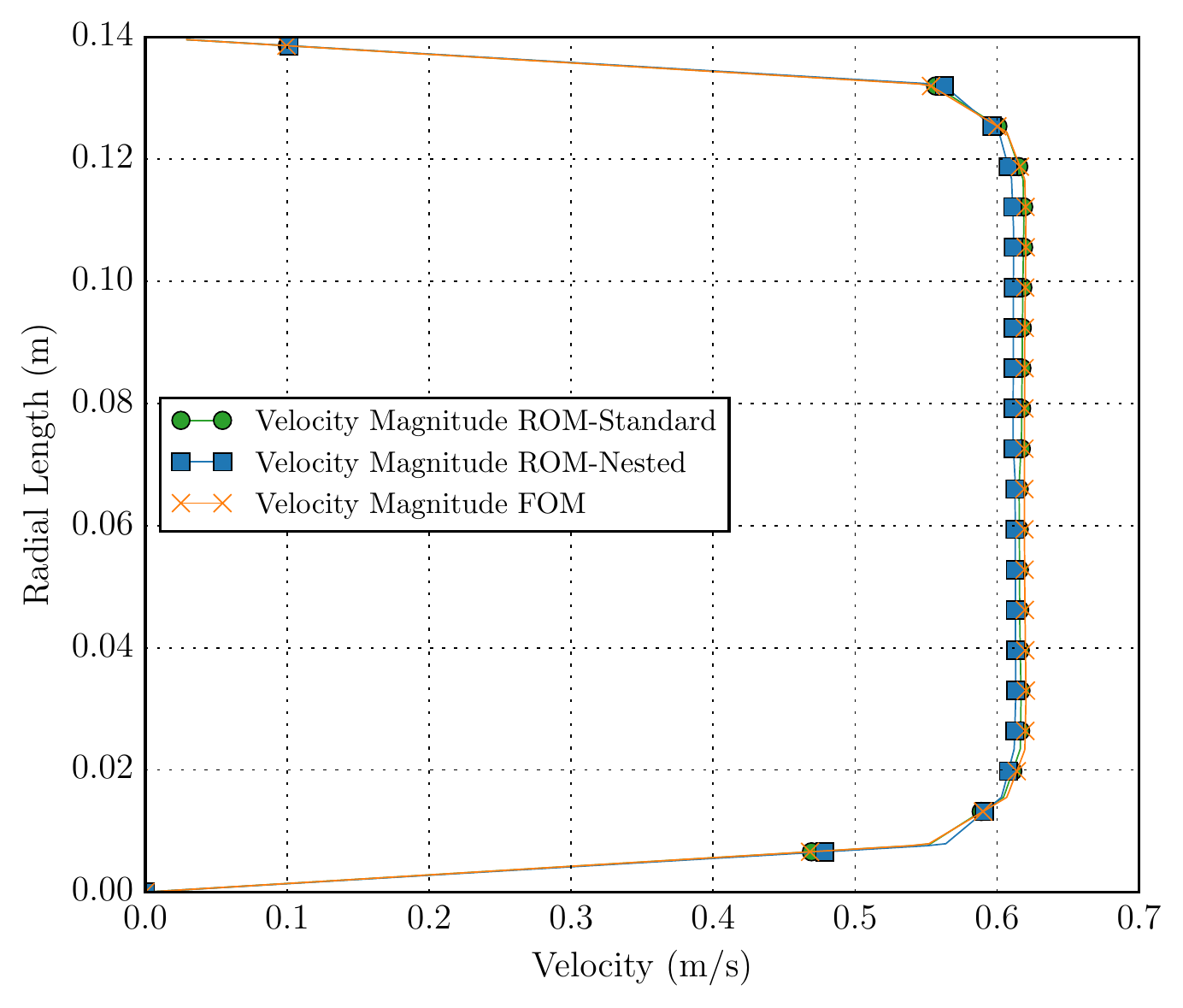}
\includegraphics[width=0.30\textwidth]{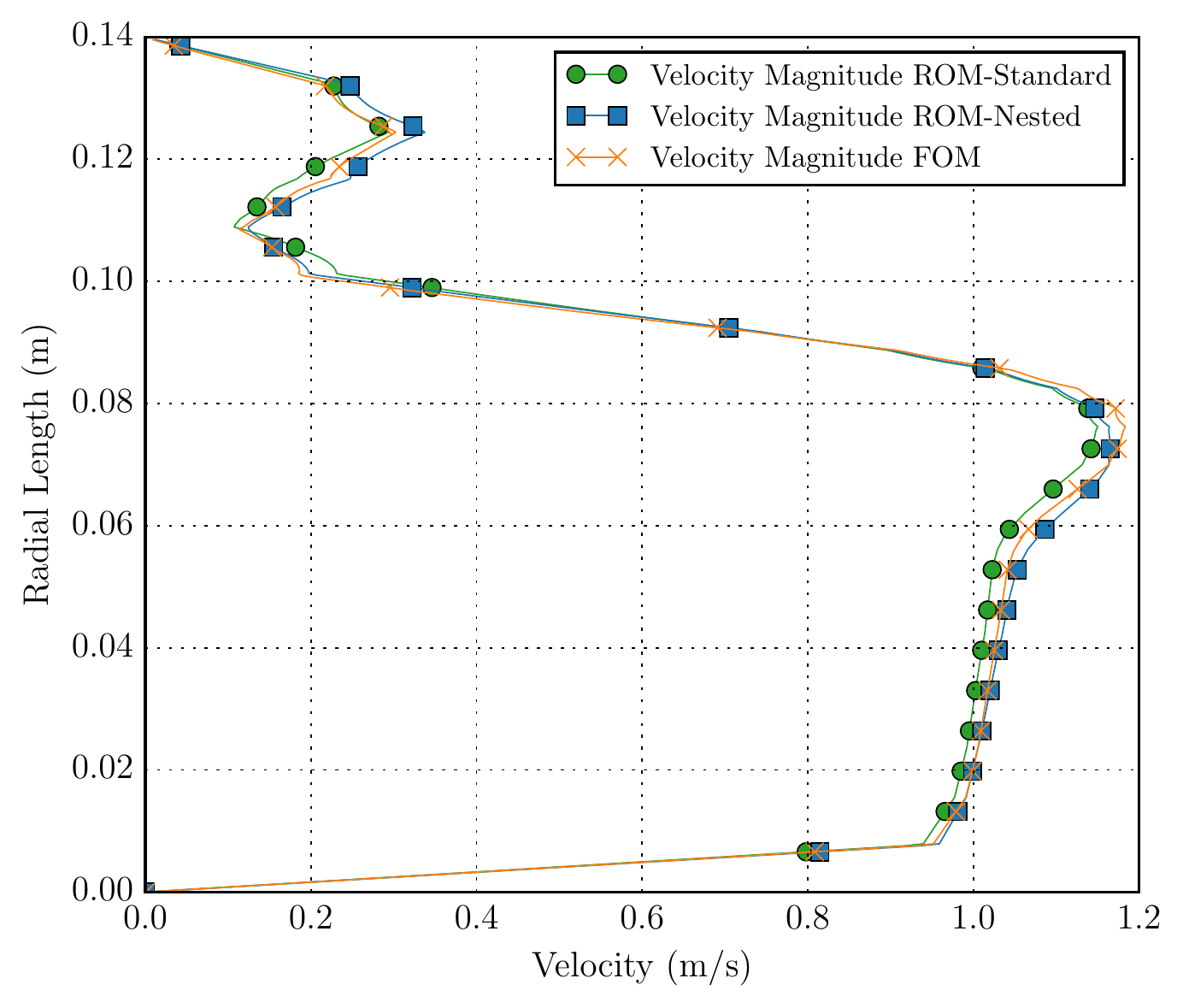}
\includegraphics[width=0.30\textwidth]{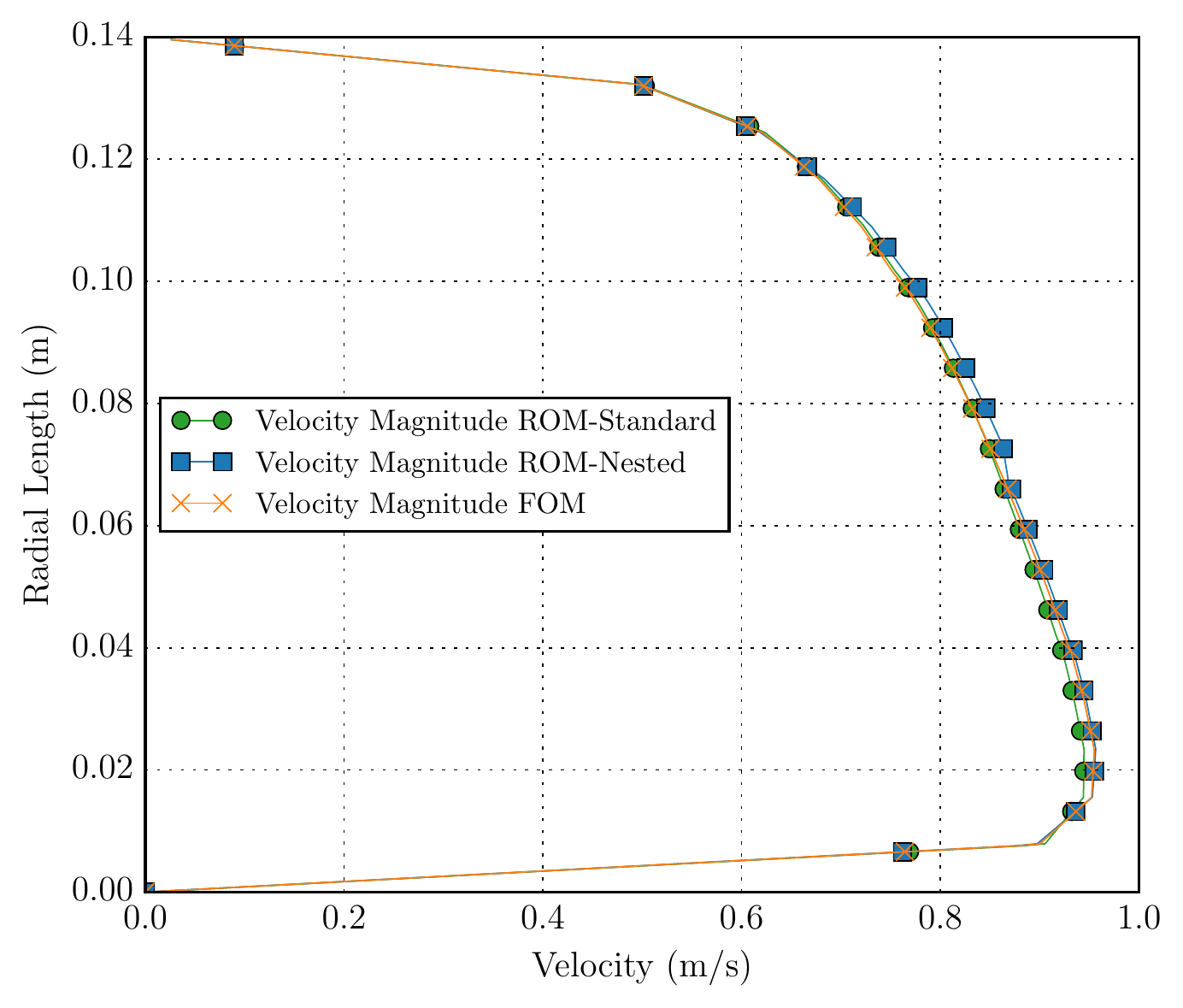}
\includegraphics[width=0.30\textwidth]{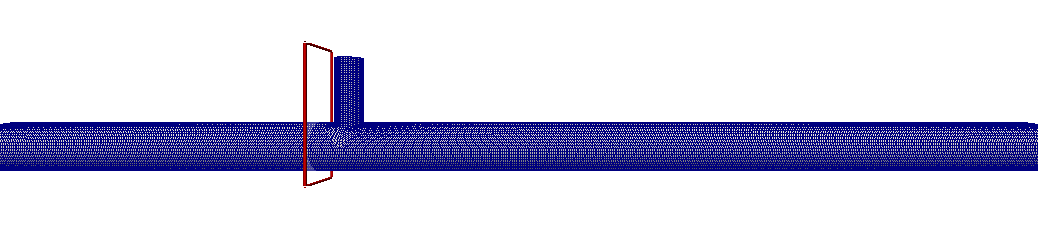}
\includegraphics[width=0.30\textwidth]{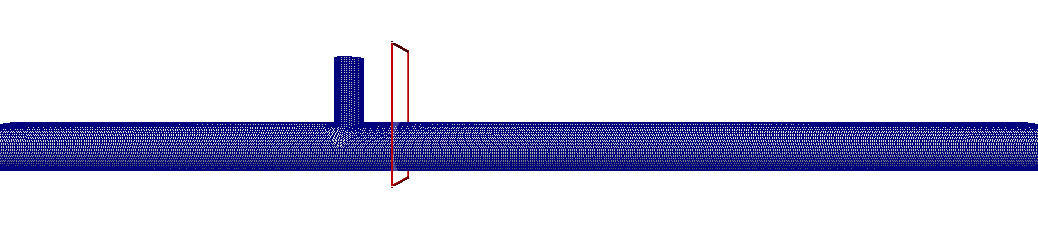}
\includegraphics[width=0.30\textwidth]{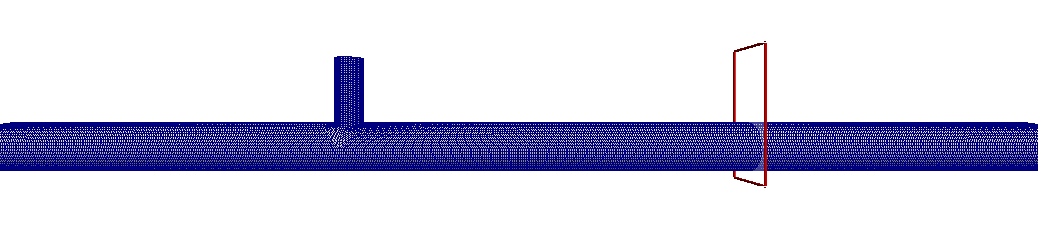}
\ 
\end{minipage}     
\caption{Comparison of the radial velocity between the FOM, ROM-Global and ROM-Nested for three different locations, before the mixing, in the mixing region and near the outlet for $t=3 \si{s}$. The second row depicts the locations where the sampling takes place.}\label{fig:radial_plots}
\end{figure*}

A visualisation of the instantaneous fields is shown in figures (\ref{fig:comparison_case2_vel}, \ref{fig:comparison_case2_temp}, \ref{fig:comparison_case2_press}) and (\ref{fig:comparison_case2_nut}) for time instances $t=0.5 \si{s}, 1.5\si{s}$ and $3 \si{s}$, where both methods show similar qualitative performance. The relative difference between the FOM and ROM fields is visualised in figures (\ref{fig:comparison_case2_diff}) and (\ref{fig:comparison_case2_diff_nut}). The biggest difference between the FOM and ROM fields is found in the mixing area of the pipe. This behaviour is expected, as in the mixing region, the flow is complex and highly transient. This behaviour is also apparent in figure (\ref{fig:radial_plots}), where the radial velocity is plotted against the arc length of the pipe in three regions: before the mixing, in the mixing and after the mixing region, close to the outlet. The mixing region is where the radial velocity diverges the most.

\begin{figure*}[!tbp] 
  \centering 
  \begin{minipage}[b]{0.40\textwidth}
    \includegraphics[width=\textwidth]{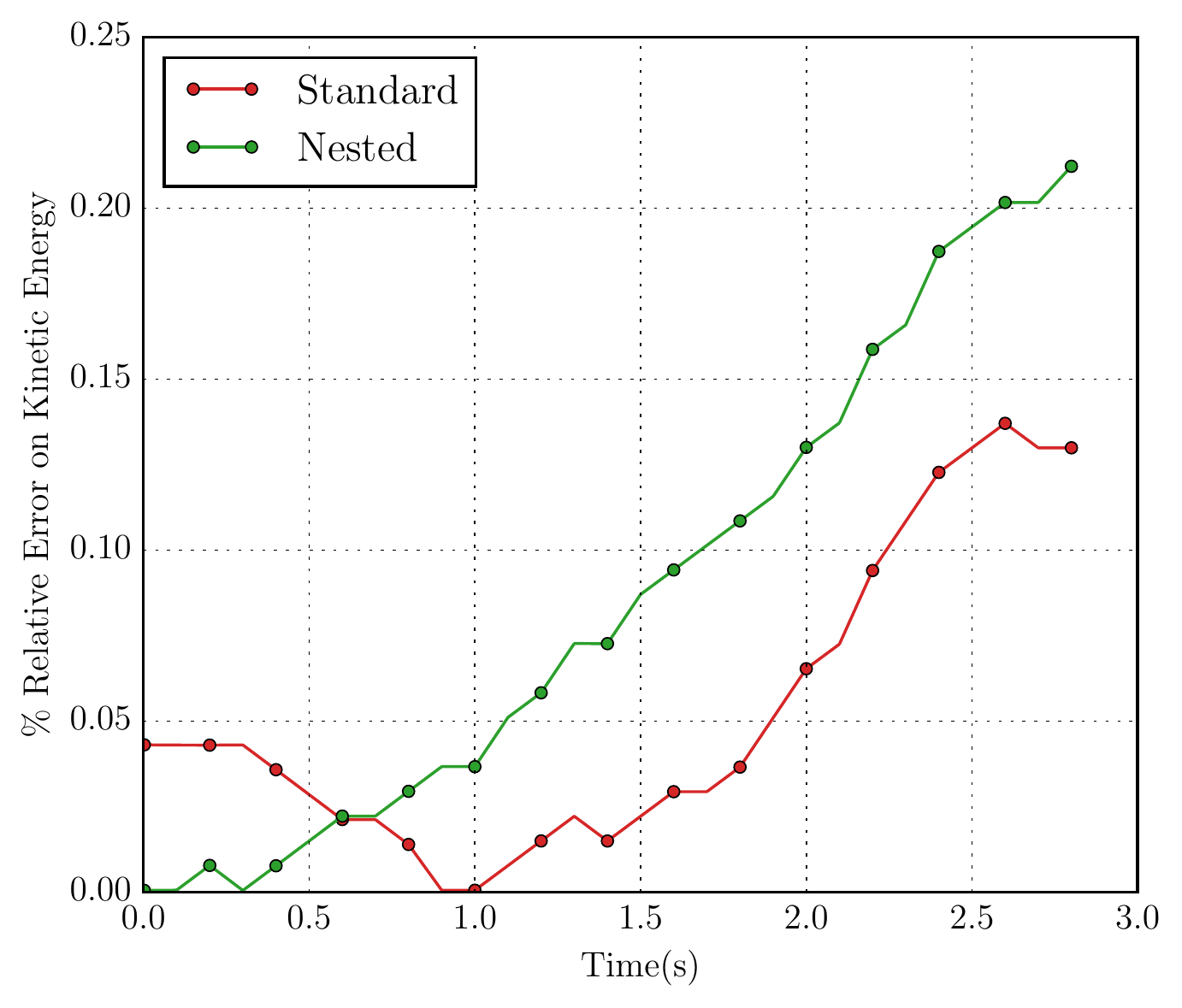}
    \caption{ Relative error on total energy (kinetic and thermal) between the FOM and the ROM for standard and nested POD methods.}\label{fig:totenergy_error}
  \end{minipage} 
\end{figure*}

To further assess the performance of the ROM against the FOM, the relative error of the total energy (kinetic and thermal) is plotted in figure (\ref{fig:totenergy_error}), where it shows a small (less than 0.25$\%$) relative error. Overall, both the ROM derived using standard and the ROM derived using nested POD methods are performing well throughout the simulation.

\begin{table*}[!tbp]
\centering
\begin{tabular}{ l | c | c | c | c}
\hline\hline
& $\Gamma_{m}$ & $\Gamma_{b}$ & $\Gamma_{w}$ & $\Gamma_{o}$  \\[0.5ex] 
$\bm{u}$ & $(0.59,0,0)$ & $(0, 0, -0.77)$ & $\nabla\cdot\bm{u} = 0$ & $\nabla \bm{u} \cdot \bm{n}$ = 0\\ \hline
$p$ & $\nabla p \cdot \bm{n} = 0$ & $\nabla p \cdot \bm{n} = 0$ & $\nabla p \cdot \bm{n} = 0$ & $0$\\ \hline 
$\theta$ & $292.15$ & $309.5$ & $\nabla \theta \cdot \bm{n} = 0$ &$\nabla \theta \cdot \bm{n} = 0$\\ 
\hline
\end{tabular}
\caption{Summary of boundary conditions where $\Gamma_{m}$ is the main pipe inlet, $\Gamma_{b}$ is the branch pipe and $\Gamma_{o}$ is the outlet.}\label{fig:boundary_con}
\end{table*}

\begin{table*}[!tbp]   
\centering
\begin{tabular}{ l | c | c  }
\hline\hline  
 & FOM & ROM   \\ [0.5ex]\hline
CPU Time(s) & $13782.1$ & $11.02$   \\  
\hline
\end{tabular}
\caption{Computational time for the full order (running on a single processor) and reduced order models.}\label{tab:cpu}       
\end{table*}

\begin{table*}[!tbp]
\begin{center}
\begin{tabular}{ l | c | c | c | c }
\hline\hline
Nested & $\bm{u}$ & $\theta$ & $p$ & $\nu_t$  \\ [0.5ex]\hline
Minimum & $1.410$ & $0.140$ & $3.097$ & $1.472$  \\ 
Maximum & $4.856$ & $0.283$ & $7.299$ & $9.826$   \\  
Average & $2.252$ & $0.187$ & $5.584$ & $6.002$  \\ 
\hline
\end{tabular}
\caption{$\%$ Relative $\epsilon_{L^2}(t)$ error for velocity, temperature and pressure and eddy viscosity fields for the nested POD method.}\label{tab:errorstatnested}       
\end{center}
\end{table*}

\begin{table*}[!tbp] 
\begin{center} 
\centering
\begin{tabular}{ l | c | c | c | c }
\hline\hline
Standard & $\bm{u}$ & $\theta$ & $p$  & $\nu_t$  \\ [0.5ex]
\hline
Minimum  & $1.411$ & $0.111$ & $3.324$ & $1.248$ \\ 
Maximum  & $2.012$ & $0.267$ & $7.996$ & $9.293$  \\ 
Average  & $1.642$ & $0.169$ & $5.563$ & $5.297$  \\ 
\hline 
\end{tabular} 
\caption{$\%$ Relative $\epsilon_{L^2}(t)$ error for velocity, temperature and pressure and eddy viscosity fields for two the standard POD method.}\label{tab:errorstat} 
\end{center}     
\end{table*}

\setlength\intextsep{0pt}

\section{Conclusion}\label{sec:concl}
In this work a hybrid reduced order model for modelling turbulent heat transfer problems has been studied. The hybrid method consists of POD-Galerkin approach for the velocity, temperature and pressure fields and PODI with RBF interpolation for the eddy viscosity field. Furthermore, two variations of POD, standard and nested, have been studied and compared. The proposed method is tested on a T-junction pipe where turbulent thermal mixing takes place. According to the results, the ROM constructed using standard or nested POD is capable of reproducing the FOM results with good accuracy in both cases, while the nested POD method requires less numerical effort. The relative error between the FOM and ROM, which peaks in the mixing region, is within acceptable levels, from engineering point of view. A speed-up factor of approximately 1259 has been obtained, meaning that the ROM is three orders of magnitude faster than the FOM. 

As future suggestions, considering the two way coupling between velocity and temperature equations would result to a ROM which could model more realistic industrial flows \cite{kstar,VERGARI2020103071}. As turbulent flow is highly complex in nature, to enhance the results and reduce the relative error, a method where local ROMs are constructed for each value of the sampling parameter and then interpolated would probably result in a more accurate ROM and would be of great interest.

\section*{Acknowledgements}     
We acknowledge the financial support of Rolls-Royce, EPSRC, the European Research Council Executive Agency by means of the H2020 ERC Consolidator Grant project AROMA-CFD ``Advanced Reduced Order Methods  with  Applications  in  Computational  Fluid  Dynamics'' - GA  681447, (PI: Prof. G. Rozza) and INdAM-GNCS 2019.

\section*{Appendix A. List of abbreviations and symbols}
\printnomenclature
\Urlmuskip=0mu plus 1mu\relax
\bibliographystyle{amsplain}
\bibliography{bib/bibfile} 
\end{multicols}
\end{document}